\documentclass[letterpaper, 11pt]{amsart}
\pdfoutput=1

\usepackage{amssymb, amsmath}
\usepackage[english]{babel}
\usepackage{ifthen}
\usepackage{ifpdf}
\usepackage{graphicx}
\usepackage{url}

\usepackage{amscd}
\usepackage{cite}

\usepackage[pdftex, plainpages = false, pdfstartview=FitH,
    pdfpagelabels,  pdfpagelayout = SinglePage, 
    bookmarks = true, bookmarksopen = true, 
    bookmarksnumbered = true, linktocpage,
    colorlinks = true, 
    linkcolor = blue, urlcolor  = blue, 
    citecolor = red,  anchorcolor = green, 
    hyperindex = true, hyperfigures]{hyperref} 

\DeclareGraphicsExtensions{.png, .jpg, .pdf}
\graphicspath{{img/}{img/}{img/}}

\title[Submanifold arrangements]%
			{Arrangements of submanifolds and the tangent bundle complement}
\author[P\, Deshpande ]{Priyavrat Deshpande}
\address{Chennai Mathematical Institute \\ Chennai\\ India}
\email{pdeshpande@cmi.ac.in} 


\def\ds{\displaystyle}

\newcommand{\ol}{\overline}

\newcommand{\maxo}{\overline{\omega}}
\newcommand{\mino}{\underline{\omega}}

\newcommand{\beq} {\begin{eqnarray}}
\newcommand{\eeq}{\end{eqnarray}}

\def\@begintheorem#1#2{\par\bgroup{\sc #1\ #2. }\it \ignorespaces}
\def\@opargbegintheorem#1#2#3{\par\bgroup{\sc #1\ #2\ (#3).}\it \ignorespaces}
\def\@endtheorem{\egroup}

\numberwithin{equation}{section}

\theoremstyle{plain}
\newtheorem{theorem}{Theorem}[section]
\newtheorem{lemma}[theorem]{Lemma}
\newtheorem{cor}[theorem]{Corollary}
\newtheorem{prop}[theorem]{Proposition}

\theoremstyle{definition}
\newtheorem{defn}[theorem]{Definition}
\newtheorem{ex}[theorem]{Example}
\newtheorem{xca}[theorem]{Exercise}
\newtheorem{conj}[theorem]{Conjecture}

\theoremstyle{remark}
\newtheorem{rem}[theorem]{Remark}

\newcommand{\bt}[1]{\begin{theorem}\label{#1}}
\newcommand{\bc}[1]{\begin{cor}\label{#1}}
\newcommand{\bl}[1]{\begin{lemma}\label{#1}}
\newcommand{\bp}[1]{\begin{prop}\label{#1}}
\newcommand{\be}[1]{\begin{ex}\label{#1}}
\newcommand{\bd}[1]{\begin{defn}\label{#1}}
\newcommand{\br}[1]{\begin{rem}\label{#1}}
\newcommand{\bx}[1]{\begin{xca}\label{#1}}
\newcommand{\bcon}[1]{\begin{conj}\label{#1}}

\newcommand{\et}{\end{theorem}}
\newcommand{\ec}{\end{cor}}
\newcommand{\el}{\end{lemma}}
\newcommand{\ep}{\end{prop}}
\newcommand{\ee}{\end{ex}}
\newcommand{\ed}{\end{defn}}
\newcommand{\exc}{\end{xca}}
\newcommand{\er}{\end{rem}}
\newcommand{\econ}{\end{conj}}

\newcommand{\bpr}{\begin{proof}}
\newcommand{\epr}{\end{proof}}

\def\A  {\mathcal{A}}
\def\FA {\mathcal{F(A)}}
\def\Ch {\mathcal{C(\A)}}

\def \F {\mathcal{F}}
\def \E {\mathcal{E}}

\def \D {\mathcal{D}}

\def \s {\mathcal{S}}
\def \L {\mathcal{L}}

\def\R{\mathbb{R}}

\def\C{\mathbb{C}}
\def\Z{\mathbb{Z}}

\begin{document}

\begin{abstract}
Drawing parallels with hyperplane arrangements, we develop the theory of arrangements of submanifolds. Given a smooth, finite dimensional, real manifold $X$ we consider a finite collection $\mathcal{A}$ of locally flat, codimension-$1$ submanifolds that intersect like hyperplanes. To such a collection we associate two combinatorial objects: the face category and the intersection poset. We also associate a topological space to the arrangement called the tangent bundle complement. It is the complement of union of tangent bundles of these submanifolds inside the tangent bundle of the ambient manifold. Our aim is to investigate the relationship between the combinatorics of the arrangement and the topology of the complement. In particular we show that the tangent bundle complement has the homotopy type of a finite cell complex. We generalize the classical theorem of Salvetti for hyperplane arrangements and show that this particular cell complex is completely determined by the face category. 
\end{abstract}
 \maketitle

\section*{Introduction}

An arrangement of hyperplanes is a finite set $\A$ of hyperplanes in $\R^l$. These hyperplanes and their intersections induce a stratification of the ambient space. These strata form a poset when ordered by topological inclusion known as the \textit{face poset}. The set of all possible intersections also forms a poset known as the \textit{intersection poset} (usually ordered by reverse inclusion). These two posets contain combinatorial information about the arrangement. The topological spaces associated with an arrangement $\A$ are the \emph{real complement} $\Ch$ and the \emph{complexified complement} $M(\A)$. The real complement is the complement of the union of hyperplanes in $\R^l$, whereas the complexified complement is the complement of the union of the complexified hyperplanes in $\C^l$. One of the aspects of the theory of arrangements is to understand the interaction between the combinatorial data of an arrangement and the topology of these complements. For 
example, one would like to comprehend to what extent the combinatorial data of an arrangement control the topological invariants, such as (co)homology or homotopy groups etc., of these complements. \par

In this paper we introduce a generalization of real hyperplane arrangements which we call the \emph{arrangements of submanifolds of codimension-$1$}. We consider situations in which finitely many submanifolds of a given manifold intersect in a way that the local information is same as that of a hyperplane arrangement but the global picture is different. Intuitively, it means that for every point of the manifold there exists a coordinate neighborhood homeomorphic to an arrangement of real hyperplanes.  We also introduce an analogue of the complexified complement in this new setting and call it \emph{the tangent bundle complement}. This paper is an attempt to answer the following question: \emph{how does the combinatorics of the intersections of these submanifolds help determine the topology of the tangent bundle complement?}\par

The topological study of the complement of a complexified hyperplane arrangement can be traced back to the work of Arnold on braid groups and configuration spaces. His results led Brieskorn to conjecture that the (complexified) complement of a real reflection arrangement is a $K(\pi, 1)$ space (i.e., its universal cover is contractible). In 1973, Deligne settled this conjecture in \cite{deli72}. His first step was to parametrize the universal cover of the complement by sequences of adjacent chambers. He showed that these galleries can be expressed in a particular normal form. This normal form was the main ingredient in proving that the constructed universal cover was contractible.\par

The combinatorial nature of Deligne's arguments led to the search for {combinatorial models for the complement}, i.e., cell complexes built using the combinatorial data of an arrangement that are homotopy equivalent to the complement. One such model, the \emph{Salvetti complex}, was introduced by Salvetti in \cite{sal1}. He showed that the incidence relations in the underlying face poset are sufficient to build this CW-complex. His construction was generalized to complex arrangements by Bj\"orner and Ziegler \cite{bz92} where the authors also showed that the face poset determines the homeomorphism type of the complement. Combinatorial models for the covering spaces have also been constructed, see for example Delucchi \cite{del1} and Paris \cite{pa2}. \par

The crux of this paper is the generalization of Salvetti's result. We prove that the incidence relations among the faces of a submanifold arrangement determine a cell complex that has the homotopy type of the tangent bundle complement (Theorem \ref{thm31}). We should also note that in the initial version of this paper we only considered the so-called combinatorially special arrangements (Section \ref{sec:csarr}). Using a categorical formulation of the Salvetti complex by d'Antonio and Delucchi in \cite{dantonio_delucchi_2011} and also using the language of cellularly stratified spaces developed by Tamaki in \cite{tamaki01, tamaki02} we have extended our proof to the current setting.\par

The paper is organized as follows. We introduce the necessary background material in Section \ref{bg}. In Section \ref{ch3def} we define the arrangement of submanifolds and provide some examples.  We introduce the tangent bundle complement in Section \ref{TBComplement} and prove that it has the homotopy type of a finite cell complex. We show that this cell complex is determined by the combinatorics of the incidence relations in the stratification induced by submanifold intersections. In Section \ref{sec:csarr} we focus on a particularly nice class of submanifold arrangements for which the combinatorial properties are similar to those of a zonotope. In particular, in Section \ref{mhcSection}, we use the theory of metrical-hemisphere complexes to describe the combinatorics. We end the paper in Section \ref{endremarks} by laying out future research and potential applications.\par

\subsection*{Acknowledgements} 
This paper is a part of the author's doctoral thesis \cite{deshpande_thesis11}. The author would like to thank his supervisor Graham Denham for his support. This particular generalization of hyperplane arrangements and the complexified complement is due to Tduesz Januskeweisz and Richard Scott. It is a pleasure to acknowledge the discussions with Tduesz Januskeweisz, without which this work would not have been possible. Sincere thanks to Mario Salvetti for being the external examiner of the thesis and also for the warm hospitality in Pisa.

\section{Background}\label{bg}
We start by briefly reviewing some basic facts about hyperplane arrangements. 

\subsection{Hyperplane Arrangements} \label{sec:prelim}

Hyperplane arrangements arise naturally in geometric, algebraic and combinatorial instances. They occur in various settings such as finite dimensional projective or affine spaces defined over a field of any characteristic. Here we formally define hyperplane arrangements and the associated combinatorial and topological data in a setting that is most relevant to our work. 

\bd{def1} A real \emph{arrangement of hyperplanes} is a collection $\A = \{H_1,\dots,H_k\}$ of finitely many \emph{hyperplanes} in $\R^l$, $l\geq 1$. \ed

The \textit{rank of an arrangement} is the largest dimension of the subspace spanned by the normals to the hyperplanes in $\A$. We assume that all our hyperplane arrangements are \textit{essential}, i.e., their rank is equal to the dimension of the ambient vector space. An arrangement is called \textit{central} if the intersection of all the hyperplanes in $\A$ is non-empty. For a subset $X$ of $\R^l$, the \emph{restriction} of $\A$ to $X$ is  the sub-arrangement $\A_X := \{H\in\A ~|~ X\subseteq H\}$.  \par

There are two posets associated with $\A$, namely, the intersection poset and the face poset which contain important combinatorial information about the arrangement. 

\bd{def2} The \emph{intersection poset} $L(\A)$ of $\A$ is defined as the set of all intersections of hyperplanes ordered by reverse inclusion. \ed

$L(\A)$ is a ranked poset with the rank of an element being the codimension of the corresponding intersection. In general it is a (meet) semilattice; it is a lattice if and only if the arrangement is central. \par 

The hyperplanes of $\A$ induce a stratification of $\R^l$ such that each stratum is an open polyhedron; these strata are called the \emph{faces} of $\A$.

\bd{def3} The \emph{face poset} $\FA$ of $\A$ is the set of all faces ordered by topological inclusion: $F\leq G$ if and only if $F\subseteq\overline{G}$. \ed

The Codimension-$0$ faces are called \emph{chambers}. The set of all chambers is denoted by $\Ch$. A chamber is \textit{bounded} if it is a bounded subset of $\R^l$. Two chambers $C$ and $D$ are \textit{adjacent} if their closures have a non-empty intersection containing a codimension-$1$ face.
By a \textit{complexified hyperplane arrangement} we mean that an arrangement of hyperplanes in $\C^l$ for which the defining equation of each hyperplane is real. Hence to every hyperplane arrangement in $\R^l$ there corresponds an arrangement of hyperplanes in $\C^l$. An important topological space associated with a real hyperplane arrangement is the following 

\bd{def4} Let $\A$ denote a hyperplane arrangement in $\R^l$. The \emph{complexified complement} of such an arrangement is denoted by $M(\A)$ and defined as 
\[M(\A) := \C^l \setminus \ds (\bigcup_{H\in\A} H_{\C}) \] 
where $H_{\C}$ is the hyperplane in $\C^l$ with the same defining equation as $H$.\ed 

Note that since $M(\A)$ is of real codimension $2$ in $\C^l$, it is connected. In fact it is an open submanifold of $\C^l$ with the homotopy type of a finite dimensional CW complex \cite[Section 5.1]{orlik92}. 

\subsection{The Salvetti Complex} \label{sec:TheSalvettiComplex}

For a hyperplane arrangement $\A$ in $\R^l$ there is a construction of a regular CW-complex, introduced by Salvetti \cite{sal1}. This complex has the same homotopy type as that of the complexified complement. The construction relies on the incidence relations among the faces. \par 

We explain the construction of the regular $l$-complex, called the Salvetti complex and denoted by $Sal(\A)$, by first describing its cells. The $k$-cells, for $0\leq k\leq l$, of $Sal(\A)$ are in one-to-one correspondence with the pairs $[F, C]$, where $F$ is a codimension-$k$ face of $\A$ and $C$ is a chamber whose closure contains $F$.\par 
Since $Sal(\A)$ is regular all the attaching maps are homeomorphisms. Hence it is enough to specify the boundary of each cell. A cell labelled $[F_1, C_1]$ is contained in the boundary of another cell labelled $[F_2, C_2]$ if and only if $F_1 \leq F_2$ in $\FA$ and $C_1, C_2$ are contained in the same chamber of $\A_{F_1}$. 

\begin{theorem}[Salvetti \cite{sal1}]\label{thm0} Let $\A$ be an arrangement of real hyperplanes and $M(\A)$ be the complement of its complexification inside $\C^l$. Then there is an embedding of $Sal(\A)$ into $M(\A)$ moreover there is a natural map in the other direction which is a deformation retraction.\et

\subsection{Topological Combinatorics}\label{tposet} 
We briefly review order complexes, the Nerve lemma and cell complexes. We assume the reader's familiarity with the basic notions related to abstract simplicial complexes and posets. The main reference is Kozlov's book \cite{dk1}. \par 

A \emph{standard $n$-simplex} is defined as the convex hull of standard unit vectors in $\R^{n+1}$. The following construction will be used as a path to move from combinatorics to topology. 

\bd{defl11} Given a finite abstract simplicial complex $\Delta$, its \emph{geometric realization} is the topological space obtained by taking the union of standard $k$-simplex in $\R^{|V|}$, for all simplices of dimension $k$. Any topological space that is homeomorphic to the realization of $\Delta$ is called the \emph{geometric realization} of $\Delta$, and is denoted by $|\Delta|$. \ed  

A topological statement about an abstract simplicial complex is in fact a statement about its geometric realization. For the sake of simplicity, in this paper, we will not differentiate between an abstract simplicial complex and its geometric realization. We will let the context decide. A simplicial map between two simplicial complexes induces a continuous map between the corresponding geometric realizations.\par
 
\bd{defoc}The \emph{order complex} $\Delta(P)$ of a poset $P$ is the abstract simplicial complex on vertex set $P$ whose $k$-faces are the $k$-chains in $P$. \ed

We will not differentiate between the order complex and its geometric realization. A poset map between two posets induces a simplicial map between the corresponding order complexes. Given a simplicial complex $\Delta$, the set of its faces ordered by inclusion forms a poset which is called the \emph{face poset}. The order complex of this face poset is the \emph{(first) barycentric subdivision} denoted by $sd(\Delta)$. The spaces $|\Delta|$ and $|sd(\Delta)|$ are homeomorphic. Thus from a topological viewpoint simplicial complexes and posets can be considered equivalent notions. \par
It is often the case that given some combinatorial data one would like to synthetically construct a topological space consistent with the data. The geometric simplicial complex is one such example. Here we consider the class of \emph{regular cell complexes}. A subset $e$ of a topological space $X$ is called a \emph{closed} (\emph{open}) $k$-\emph{cell} if it is homeomorphic to (the interior of) the standard $k$-ball  in $\R^k$.\par

\bd{def5ch1s1}
A \emph{regular cell complex} $(X, \mathcal{C})$ is a pair consisting of a Hausdorff space $X$ and a finite collection $\mathcal{C}$ of open cells in $X$ such that
\begin{enumerate}
\item $X = \bigcup_{e\in \mathcal{C}} e$,
\item the boundary $\overline{e}\setminus e$ of each cell is a union of some members of $\mathcal{C}$.
\end{enumerate} \ed

Given a (regular) CW-complex $X$, the set consisting of its cells ordered by inclusion forms a poset. This poset is the \textit{face poset} and is denoted by $\F(X)$. We state a result concerning regular cell complexes which is important from the combinatorial viewpoint. 

\bt{thm2ch1s3}
Let $(X, \mathcal{C})$ be a regular cell complex and $\F(X)$ be its face poset. Then \[\Delta(\F(X)) \cong X. \]
Furthermore, this homeomorphism can be chosen so that it restricts to a homeomorphism between $\overline{e}$ and $\Delta(\F_{\leq e})$, for all $e\in\mathcal{C}$. \et

However, we will need a more general setting which we now define.

\bd{defn1s1} A small category is called \textit{acyclic} if only identity morphisms have inverses, and any morphism from an object to itself is an identity.\ed

A poset $(P, \leq)$ can be regarded as an acyclic category as follows. The objects of this category are elements of $P$. There is a unique morphism from $x$ to $y$ if $x\leq y$. In fact, one can think of acyclic categories as generalization of posets in the following sense. Given an acyclic category $C$ there is a partial order $\prec$ on its objects given by $C(x, y)\neq \emptyset \Rightarrow x \prec y$. $(\mathrm{Ob}(C), \prec)$ is called the \textit{underlying poset} of $C$. If $C$ is an acyclic category with at most $1$ morphism between any two objects then $C$ and its underlying posets are equivalent as acyclic categories. This is precisely what it will mean when we say that an acyclic category is (equivalent to) a poset.\par 

As stated earlier the face poset of a simiplicial complex captures the incidence relations between simplices. Similarly an acyclic category encodes incidence relations among cells of more general type of cell complexes. Before we introduce those spaces let us look at a generalization of the order complex construction. We refer the reader to Kozlov's book \cite[Definition 10.4]{dk1} for the precise technical definition.

\bd{defnrevac} The nerve (or the geometric realization) of an acyclic category $C$, denoted by $\Delta(C)$, is space constructed as follows: 
\begin{itemize}
 \item The ordered vertex set of $\Delta(C)$ is the object set of $C$.
 \item For $k\geq 1$, there is a set of $k$-simplices denoted by $S_k(\Delta(C))$. This set contains all the composable morphism chains consisting of $k$ nonidentity morphisms.  
\end{itemize}
A $(k-1)$-simplex is glued to a face of a $k$-simplex by the canonical linear homeomorphism between them that preserves the ordering of the vertices. 
\ed

The geometric realization of an acyclic category is a \textit{regular trisp} (\textit{delta complex}). Recall that a delta complex is a quotient space of a collection of disjoint simplices obtained by identifying certain of their faces via canonical linear homeomorphisms that preserve the ordering of vertices. In the case of regular trisps there are no self identifications in the boundary of a simplex. See \cite[Section 2.4]{dk1} for details.\par 

\bd{defn2s1} For an arbitrary acyclic category $C$, let $sd C$ denote the poset whose minimal elements are objects of $C$, and whose other elements are all composable morphism chains consisting of nonidentity morphisms of $C$. The elementary order relations are given by composing morphisms in the chain, and by removing the first or the last morphism. The poset $sd C$ is called the \textit{barycentric subdivision} of $C$.\ed
\br{remacyclic}The barycentric subdivision of an acyclic category can also be defined as the face poset of its nerve. Hence there is a regular CW complex homeomorphic to the geometric realization of $C$.\er 

The \emph{Nerve Lemma} often helps simplify a given topological space for combinatorial applications. Recall that an \emph{open covering} of a topological space is a collection of open subsets of the space such that their union is the space itself. The \emph{nerve of an open covering} is a simplicial complex whose vertices correspond to the open sets and a set of $k+1$ vertices spans a $k$-simplex whenever the corresponding $k+1$ open sets have a nonempty intersection. In general the nerve need not reflect the topology of the ambient space, but the following result gives a useful condition when it does. 

\begin{theorem}[Theorem 15.21\cite{dk1}]\label{thm1ch1s1} If $\mathcal{U}$ is a finite open cover of a topological space $X$ such that every non empty intersection of open sets in $\mathcal{U}$ is contractible, then $X\simeq \mathrm{nerve}(\mathcal{U})$.  \et

\subsection{Cellularly stratified spaces}
We now introduce the language of cellular stratified spaces developed by Tamaki in \cite{tamaki01, tamaki02}. A hyperplane arrangement induces a stratification of the ambient space such that each stratum is a relatively open, convex polyhedron. These strata might be unbounded regions. Consequently such a stratification is not regarded as a CW-complex structure in the usual sense. One needs to appropriately modify the definition so that the induced stratification behaves like a regular CW-complex. 

\bd{defn3s1}Let $X$ be a topological space and $P$ be a poset. A stratification of $X$ indexed by $P$ is a surjective map $\sigma\colon X\to P$ satisfying the following properties:
\begin{enumerate}
	\item For $p\in P$, $e_p := \sigma^{-1}(p)$ is connected and locally closed.
	\item For $p, q\in P,~ e_p\subseteq \ol{e_q} \iff p\leq q$.
	\item $e_p \cap \ol{e_q} \neq \emptyset \Longrightarrow e_p \subseteq \ol{e_q}$.
\end{enumerate}
The subspace $e_p$ is called as the stratum with index $p$. \ed

\br{remloclosed}A subset $A$ of a topological space $X$ is \textit{locally closed} if $A$ is open in $\ol{A}$.\er

The reader can verify that the boundary of each stratum is itself a union of strata. Such a stratification defines a decomposition of $X$. The poset $P$ is called the \textit{face poset} of $X$. A \textit{morphism} of stratified spaces is a pair $\mathbf{f} = (f, \underline{f})\colon (X, \sigma_X)\to (Y, \sigma_Y)$ consisting of a continuous map $f\colon X\to Y$ and a poset map $\underline{f}$ making the obvious diagram commute. $f$ is \textit{strict} if $f(e_p) = e_{\underline{f}(p)}$ for every $p \in \mathrm{Im} \sigma_X$. A subspace $A\subseteq X$ is a \textit{stratified subspace} if the restriction $\sigma |_{A}$ is a stratification. When the inclusion of $A$ into $X$ is a strict morphism, $(A, \sigma |_{A})$ is called a \textit{strict stratified subspace}.\par 
In general, each stratum of a stratified space could have non-trivial topology, however, we focus on those stratifications in which strata are cells.

\bd{defn4s1}
A \textit{globular $n$-cell} is a subset $D$ of $D^n$ (the unit $n$-disk) containing $\mathrm{Int} (D^n)$. We call $D\cap \partial(D^n)$
the boundary of $D$ and denote it by $\partial D$. The number $n$ is called the (\textit{globular}) \textit{dimension} of $D$.\ed

\bd{defn5s1} Let $X$ be a Hausdorff space. A cellular stratification of $X$ is a pair $(\sigma, \Phi)$ of a stratification $\sigma$ and a finite collection of continuous maps $\Phi = \{\phi_p : D_p\to \ol{e_p}~|~ p\in P\}$, called cell structures, satisfying following conditions:
\begin{enumerate}
	\item $\phi_p(D_p) = \ol{e_p}$ and $\phi_p|_{\mathrm{Int}(D^n)}\colon \mathrm{Int}(D^n)\to e_p$ is a homeomorphism.
	\item $\phi_p$ is a quotient map for every $p\in P$. 
	\item For each $n$-cell $e_p$ the boundary $\partial e_p$ contains cells of dimension $i$ for $0\leq i\leq n-1$.
\end{enumerate} 
A cellularly stratified space is a triple $(X, \sigma, \Phi)$, where $(\sigma, \Phi)$ is a cellular stratification. \ed

We assume that all our cellularly stratified spaces (CS-spaces for short) have a finite number of cells. Consequently the cell structure is CW (i.e., each cell meets only finitely many other cells and X has the weak topology determined by the union of cells). Morphisms of CS-spaces and strict cellular stratified subspaces can be defined analogously. A cell $e_p$ is said to be regular if the structure map $\phi_p$ is a homeomorphism. A cell complex is regular if all its cells are regular.

\bd{defn6s1} Let $X$ be a CS-space. $X$ is called totally normal if, for each $n$-cell $e_p$,
\begin{enumerate}
	\item there exists a structure of regular cell complex on $S^{n-1}$ containing $\partial D_p$ as a strict cellular stratified subspace of $S^{n-1}$ and
	\item for any cell $D$ in $\partial D_q$, there exists $e_q\in \partial e_p$ such that $D\cong D_q$ and $\phi_q = \phi_p|_{\ol{D}}$.
\end{enumerate} \ed

In particular, the closure of each $k$-cell contains at least one cell of dimension $i$, for every $i\leq k-1$. Consider $\mathrm{Int}(D^2)\cup \{(0, 1)\}$; its boundary does not contain a $1$-cell. This is an example of a CS-space which is not totally normal.  

\bd{remn7s1} The \textit{face category} of a cellularly stratified space $(X, \sigma, \Phi)$ is denoted by $\F(X)$. The objects of this category are the cells of $X$. For each pair $p \leq q$, define $\F(X)(e_p, e_q)$ (i.e., the set of all morphisms) to be the set of all maps $b: D_p\to D_q$ such that $\phi_q\circ b = \iota \circ \phi_p$.  \ed

Note that, in general, the set of morphisms between two objects in the face category has a non-trivial topology (see \cite[Definition 4.1]{tamaki01}). However, we assume, without loss of generality, that the set morphisms between any two objects is always finite. 

\begin{lemma}[Lemma 4.2 \cite{tamaki01}]\label{lem4t2} If $X$ is a totally normal CS-space then its face category $\F(X)$ is acyclic. Moreover if $X$ is regular then $\F(X)$ is equivalent to the dual of the face poset of $X$.\el

For CS-spaces the poset underlying its face category coincides with the (classical) face poset of a cell complex. The reason we are focusing on these spaces is the following lemma.

\bl{lemn0s1} Let $\A$ be an arrangement of hyperplanes in $\R^l$. Then the stratification induced by the hyperplanes in $\A$ defines the structure of a totally normal, cellularly stratified space. \el


As an example consider the arrangement of coordinate axes in $\R^2$. For the closed cell $e = \{(x,y) ~|~ x, y\geq 0 \}$ the subset 
$\mathrm{Int}D^2\cup \{(x,y)\in S^1~|~ x < 0 \}$ can serve as the globular cell. \par 

We state one more property of these spaces which we wish to use in order to define arrangements of submanifolds. In particular, the following result generalizes Theorem \ref{thm2ch1s3}.
 
\begin{lemma}[Corollary  4.17\cite{tamaki01}] \label{lemn1s1} For a totally normal cellularly stratified space $X$, the nerve of its face category embeds in $X$ as a strong deformation retract. \el

\br{rem_dual}
If $(X, \sigma, \Phi)$ is a totally normal CS-space then one can also define its dual $X^* := (X, \sigma^*, \Phi^*)$. The construction is such that $X^*$ is also a totally normal CS-space. There is an order-reversing poset isomorphism between the face poset of $X^*$ and that of $X$. The attaching maps in $\Phi^*$ are defined such that the face category of $X^*$ is the opposite category of $\F(X)$. \er

\subsection{Homotopy theoretic techniques}
Now we state relevant results from homotopy theory applied to combinatorial settings. The approach we want to take here is more concrete; we follow the paper of Welker, Ziegler and \u{Z}ivaljevi\'{c} \cite{wzz}, see also \cite[Chapter 15]{dk1}. All of their results are for diagrams indexed by posets. We modify the results to our setting. The crucial observation is that in most of their results one can replace posets by acyclic categories. \par 

A \emph{diagram of spaces} indexed over an acyclic category is a covariant functor $\D\colon C\to \mathrm{Top} $ from an acyclic category $C$ to the category topological spaces and continuous maps (in this case $\D$ will be called a $C$-diagram).

\bd{def362} Consider a diagram $\D\colon C\to \mathrm{Top}$, then the homotopy colimit of $\D$ is defined as 
\[\mathrm{hocolim}\D := \ds \coprod_{\sigma = v_0\to\cdots\to v_n}(\sigma\times \D(v_0))/\sim \]
where the disjoint union is taken over all simplices in $\Delta(C)$. The equivalence relation $\sim$ is defined by following condition: for $\tau_i\in \partial\sigma$, let $f_i\colon \tau_i\hookrightarrow\sigma$ be the inclusion map; then
\begin{enumerate}
 \item for $i>0$, $\tau_i\times \D(v_0)$ is identified with the subset of $\sigma\times\D(v_0)$ by the map induced by $f_i$;
 \item for $\tau_0 = v_1\to\cdots \to v_n$, we have $f_0(\alpha)\times x\sim \alpha \times\D(v_0\to v_1)(x)$ for any $\alpha\in\tau_i$ and $x\in\D(v_0)$.
\end{enumerate}
\ed
Recall that in a small category $P$ and for some $p\in \mathrm{Ob}(P)$ the comma category, denoted $P\downarrow p$, is the category whose objects are the morphisms of $P$ ending at $p$ and morphisms are commuting triangles. In its simplest form the next theorem provides sufficient conditions for a map between two posets to be a homotopy equivalence at the level of order complexes.
\begin{theorem}[Quillen's theorem A]\label{qta}
Let $f\colon P\to Q$ be a functor of acyclic categories such that $\Delta(f^{-1}(Q\downarrow q))$ is contractible $\forall q$. If $\D$ is any $Q$-diagram and $f^*\D$ the corresponding (pull-back) $P$-diagram then, 
\[\mathrm{hocolim}_{P}(f^*\D) \stackrel{\simeq}{\rightarrow}  \mathrm{hocolim}_{Q}\D. \] \et

We now look at an explicit connection between the nerve construction and the homotopy colimit construction. 
 
\bd{def37}Let $P$ be an acyclic category. A \emph{diagram of acyclic categories} is a diagram indexed over $P$ that assigns an acyclic category $Q_p$ to each $p\in P$. \ed

Given a diagram of acyclic categories $\D$ we construct a new acyclic category, called as the acyclic limit of $\D$. It is denoted by $\mathrm{Aclim} \D$ and defined as follows:
\[ \ds \mathrm{Ob}(\mathrm{Aclim} \D) := \coprod_{p\in\mathcal{P}}\{p\}\times Q_p.\]
The morphisms in this category can be described as follows: 
\[(p,q)\to(p',q')\Leftrightarrow \exists f\colon p\to p'\hbox{~such that~} \D(f)(q) = q'. \]
For a diagram of acyclic category its homotopy colimit is indeed a nerve as the following lemma suggests.
\begin{lemma}[Simplicial Model Lemma]\label{lem35} Let $\D\colon P \to \mathrm{Top}$ be a diagram of acyclic categories. Then 
\[\mathrm{hocolim} \D\simeq \Delta(\mathrm{Aclim} \D). \] \el
The above lemma is a special case of the Grothendieck construction, which is used to define the homotopy colimit of a diagram of categories.

\section{Arrangements of Submanifolds}\label{ch3def}

In this section we propose a generalization of hyperplane arrangements. In order to do this we isolate the following characteristics of a hyperplane arrangement: 
\begin{enumerate}
\item there are finitely many (translated) codimension-$1$ subspaces each of which separates $\R^l$ into two components, 
\item there is a stratification of the ambient space into open polyhedra, 
\item the face poset of this stratification has the homotopy type of $\R^l$. 
\end{enumerate}

Any reasonable generalization of hyperplane arrangements should posses these properties. Since smooth manifolds are locally Euclidean they are obvious candidates for the ambient space. In this setting we can study arrangements of codimension-$1$ submanifolds that satisfy certain nice conditions. For example, locally, we would like our submanifolds to behave like hyperplanes.\par

\subsection{Locally flat submanifolds}
We now generalize property (1) mentioned above. Throughout this paper a manifold always has empty boundary. Our focus is on the codimension-$1$ submanifolds that are \textit{embedded as a closed subset} of finite-dimensional smooth manifold. This type of submanifold behaves much like hyperplanes. We begin by exploring some of its well-known separation properties.
 
\bl{lem1news1c3}
If $X$ is a connected $l$-manifold and $N$ is a connected $(l-1)$-manifold embedded in $X$ as a closed subset, then $X\setminus N$ has either $1$ or $2$ components. If in addition $H_1(X, \Z_2) = 0$ then $X\setminus N$ has exactly two components (i.e., $N$ separates $X$). \el

\begin{proof} The lemma follows from the following exact sequence of pairs in mod $2$ homology:
\[{ H_1(X, \Z_2)\to H_1(X, X\setminus N, \Z_2)\to \tilde{H}_0(X\setminus N, \Z_2) \to \tilde{H}_0(X, \Z_2) = 0.} \qedhere\]\end{proof}

\bl{lem2news1c3} Let $X$ be a connected $l$-manifold and let $N$ be a connected, $(l-1)$-submanifold embedded in $X$ as a closed subset. Then $N$ separates $X$ if and only if the inclusion induced homomorphism $H_c^{n-1}(X, \Z_2)\to H^{n-1}_c(N, \Z_2)$ (on the cohomology with compact supports) is trivial. \el

\begin{proof}
The proof is a simple diagram chase:
\[  {\begin{CD}
 H^{n-1}_c(X, \Z_2) @> >> H^{n-1}_c(N, \Z_2) \cong \Z_2\\
    @V\cong VV @VV\cong V \\
    H_1(X, \Z_2) @> >> H_1(X, X\setminus N, \Z_2)@> >> \tilde{H}_0(X\setminus N, \Z_2) @> >> 0. 
\end{CD}}\qedhere \] 
\end{proof}

In general a codimension-$1$ submanifold need not separate $X$; there are weaker separation conditions. 

\bd{def1news1c3} A connected codimension-$1$ submanifold $N$ in $X$ is said to be \emph{two-sided} if $N$ has a neighbourhood $U_N$ such that $U_N\setminus N$ has two connected components; otherwise $N$ is said to be \emph{one-sided}. $N$ is \textit{locally two-sided }in $X$ if each $x\in N$ has arbitrarily small connected neighbourhoods $U_x$ such that $U_x\setminus (U_x\cap N)$ has two components. \ed

Note that being two-sided is a local condition. For example, a point in $S^1$ does not separate $S^1$, however, it is two-sided. The following corollary follows from the	above definition and the lemmas.

\bc{cor2news1c3} If $X$ is an $l$-manifold and $N$ is an $(l-1)$-manifold embedded in $X$ as a closed subset such that $H_1(N, \Z_2)$ is trivial then $N$ is two-sided. \ec

\bc{cor1s1c3} Every codimension-$1$ submanifold $N$ is locally two-sided in $X$. \ec


An $n$-manifold $N$ contained in an $l$-manifold $X$ is \emph{locally flat} at $x\in N$; if there exists a neighbourhood $U_x$ of $x$ in $X$ such that $(U_x, U_x\cap N)\cong (\R^l, \R^n)$. An embedding $f\colon\thinspace N\to X$ such that $f(N)\subseteq X$ is said to be \emph{locally flat at a point} $x\in N$ if $f(N)$ is locally flat at $f(x)$. Embeddings and submanifolds are \emph{locally flat} if they are locally flat at every point.\par 

It is necessary to consider the locally flat class of submanifolds otherwise one could run into pathological situations. For example, the Alexander horned sphere is a non-flat embedding of $S^2$ inside $S^3$ such that the connected components of its complement are not even simply connected see Rushing \cite[Page 65]{rushing73} for details. \par

Locally flat submanifolds need not intersect like hyperplanes. A simple example comes from the non-Pappus arrangement of $9$ pseudolines in $\R^2$. Corresponding to this arrangement there is a rank $3$ oriented matroid which realizes an arrangement of pseudocircles in $S^2\subseteq \R^3$, see Figure \ref{nonpap} below. 
 \begin{figure}[!ht]
  \begin{center}  
    \includegraphics[scale=0.2,clip]{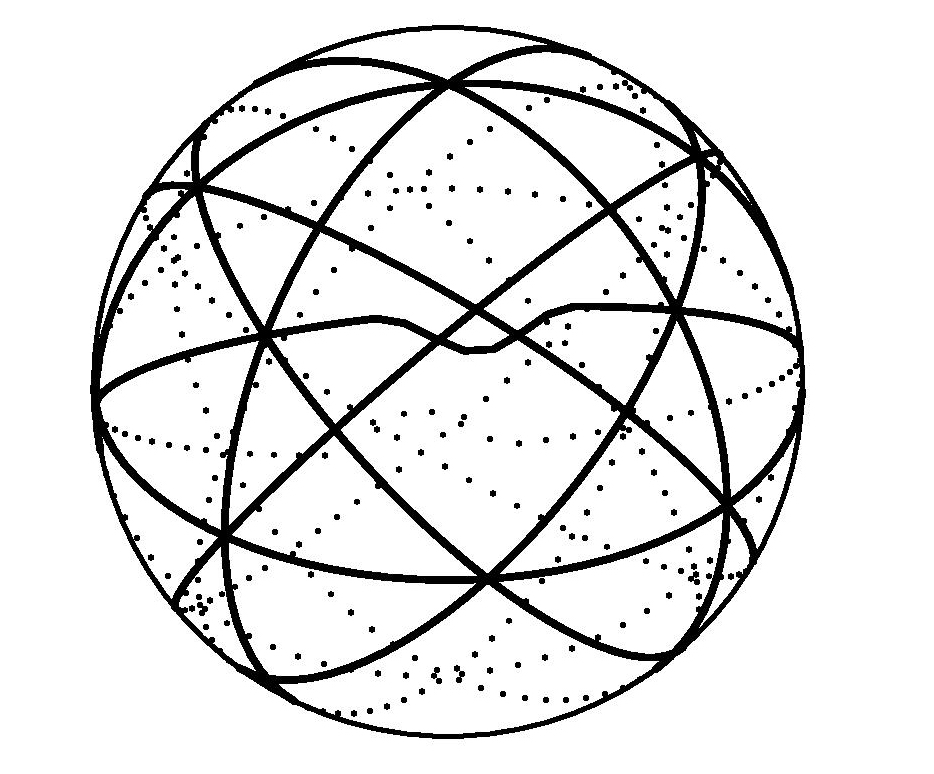} 
  \end{center}
  \caption{Non-Pappus pseudocircle arrangement.}     \label{nonpap}     
 \end{figure}

Now consider cones over these pesudocircles. The cones are pseudoplanes in $\R^3$. Each pseudoplane, individually, is homeomorphic to a $2$-dimensional subspace of $\R^3$. However there does not exist a homeomorphism of $\R^3$ which will map all of these to planes passing through a common point. For a detailed discussion about pseudo-arrangements see Bj\"orner et. al. \cite[Chapter 5]{ombook99}. \par 

We would like to avoid such situations. Hence we introduce a condition that will guarantee hyperplane-like intersections. Let $\A= \{N_1, \dots, N_k \}$ be a collection of locally flat, codimension-$1$ submanifolds of $X$. For every $x\in X$ and an open neighbourhood $V_x$ (homeomorphic to $\R^l$) of $x$ let $\A_x := \{N\cap V_x ~|~ x\in N\in \A \}$. By $\bigcup \A_x$ we mean the union of elements of $\A_x$.

\bd{ch3def0} Let $X$ be a manifold of dimension $l$. Let $\A = \{N_1,\dots,N_k\}$ be a collection of codimension-$1$, locally flat submanifolds of $X$. We say that these submanifolds have a \emph{locally flat intersection} 
if for every $x\in X$ there exists an open neighbourhood $V_x$ and a homeomorphism $\phi\colon V_x\to \R^l$ such that $(V_x, \bigcup \A_x)\cong (\R^l, \bigcup \A')$ where $\A'$ is a central hyperplane arrangement in $\R^l$ with $\phi(x)$ as a common point. \ed

\subsection{Cellular Stratifications}
Now we generalize properties (2) and (3). Let $\A = \{N_1,\dots,N_k\}$ be a collection of codimension-$1$ submanifolds of $X$ having locally flat intersection. Let $\L$ denote the set of all non-empty intersections and $\L^d$ be the subset of codimension-$d$ intersections. We have $\bigcup \L^0 = X$ and $\bigcup \L^1 = \bigcup_{i=1}^k N_i$. For each $d\geq 0$ consider the following subsets of $X$   
\begin{align*}
	\s^d(X) =  \bigcup \L^d \setminus \bigcup \L^{d+1}.
\end{align*} 

Note that each $\s^d(X)$ may be disconnected and that $X$ can be expressed as the disjoint union of these connected components. We want these components to define a `nice' stratification of $X$. For example, in case of hyperplane arrangements each stratum is an open polyhedron of appropriate dimension. We need to focus on stratifications such that the strata are cells otherwise there is no hope of generalizing property (3). For example, consider two non-intersecting longitudinal circles in the $2$-torus $S^1\times S^1$; there are two codimension-$0$ strata and two codimension-$1$ strata. The resulting face poset does not have the homotopy type of the torus. In view of Lemma \ref{lemn1s1} it is desirable to assume that the induced stratification is cellular and totally normal (Definition \ref{defn6s1}). 

\br{remcell}
Imposing this condition on the stratification implies that the resulting face category has the homotopy type of the ambient manifold (Lemma \ref{lemn1s1}). Moreover, each stratum is locally polyhedral (see \cite[Section 3.3]{tamaki01} and \cite[Lemma 2.42]{tamaki02}). Locally polyhedral means that the cell structure of each stratum is defined using a convex polyhedron and PL maps (see Lemma \ref{lemn0s1}). At first glance this might seem restrictive. For example, the stratification of $(S^1)^l$ induced by the braid arrangement is not even cellular. However, there are ways to deal with such situations. We discuss this in detail in Remark \ref{remtoric}.\er

\subsection{Definition and examples}
The desired generalization of hyperplane arrangements is the following

\bd{def31}Let $X$ be a connected, smooth, real manifold of dimension $l$. An \textbf{arrangement of submanifolds} is a finite collection 
$\A = \{N_1,\dots, N_k\}$ of codimension-$1$ smooth submanifolds in $X$ such that:
	\begin{enumerate}
		\item the $N_i$'s have locally flat intersection,
		\item the stratification induced by the intersections of $N_i$'s defines the structure of a totally normal cellularly stratified space on $X$.
	\end{enumerate}
\ed

Let us look at how we can associate combinatorial data to such an arrangement and at a few examples. 

\bd{def32}The \emph{intersection poset}, denoted by $L(\A)$, is the set of connected components of all possible intersections of $N_i$'s ordered by reverse inclusion. The rank of each element in $L(\A)$ is defined to be the codimension of the corresponding intersection.\ed

By convention $X$ is the intersection of no submanifolds, hence it is the smallest member. 

\bd{def33}The cells of the stratification will be called the \emph{faces} of the arrangement. The \textit{face category} of the arrangement is the face category of the induced stratification and is denoted by $\FA$. Codimension-$0$ faces are called \textit{chambers} and the set of all chambers is denoted by $\Ch$. \ed

\br{remn1s2} Note that the face category of an arrangement is acyclic. It is equivalent to the face poset if and only if the cell structure induced by $\A$ is regular. \er

If there is no confusion about the arrangement, we will omit $\A$ and denote the above two objects by $\F$ and $L$.\par 
The obvious examples of these submanifold arrangements are the hyperplane arrangements. Let us look at different examples. We start with the unit circle $S^1$; here the codimension-$1$ submanifolds are just points.

\be{example1s2} Consider the arrangement $\A = \{p = (1, 0) \}$ in $S^1$. The strata define the minimal cell decomposition of $S^1$; $\{p\} \cup (S^1\setminus \{p\})$. The cell structure map for the $1$-cell $\phi_1: D^1\to S^1$ can be given by $\phi_1(t) = (\cos(\pi(t+1)), \sin(\pi(t+1)))$. Moreover, there are two lifts of the cell structure map for the unique $0$-cell. Hence the face category contains two objects, the $1$-cell and the $0$-cell. Apart from the identity morphisms there are two morphisms from the $0$-cell to the $1$-cell which we denote by $\alpha_{pC}$ and $\beta_{pC}$. Figure \ref{s1pt1fig} shows the arrangement and the nerve of the corresponding face category.  
\begin{figure}[!ht]
\begin{center}  \includegraphics[scale=0.3,clip]{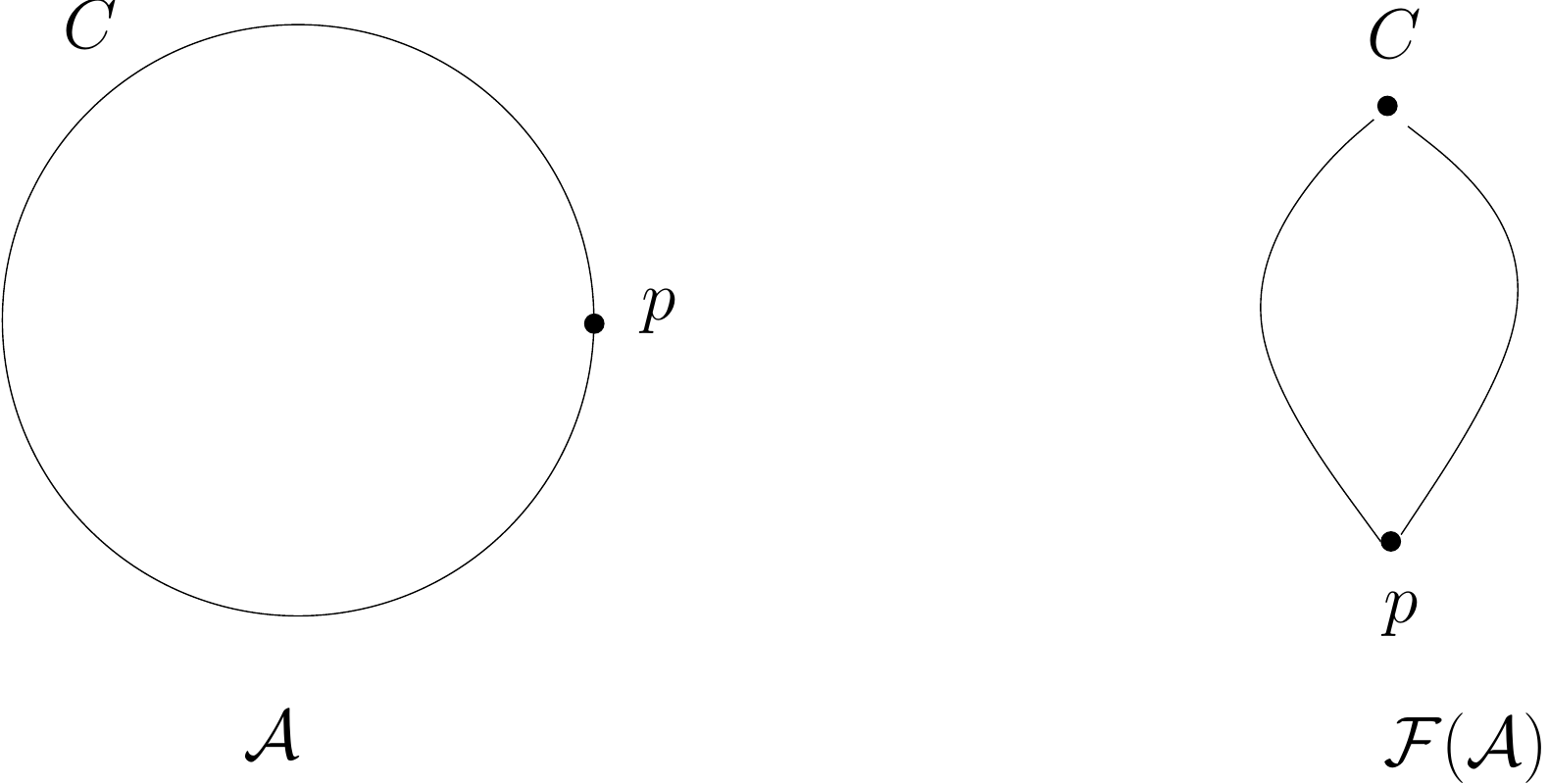} \end{center}
    \caption{Arrangement of $1$ point in a circle.}     \label{s1pt1fig}     \end{figure}
\ee

\be{ex41} Consider the arrangement $\A = \{p, q\}$ (i.e., a $0$-sphere) in $S^1$. The stratification defines a regular cell structure. Hence the face category is equivalent to the face poset. Figure \ref{figdef1} shows this arrangement and the Hasse diagrams of the face poset and the intersection poset. 
\begin{figure}[!ht]
\begin{center}  \includegraphics[scale=0.5,clip]{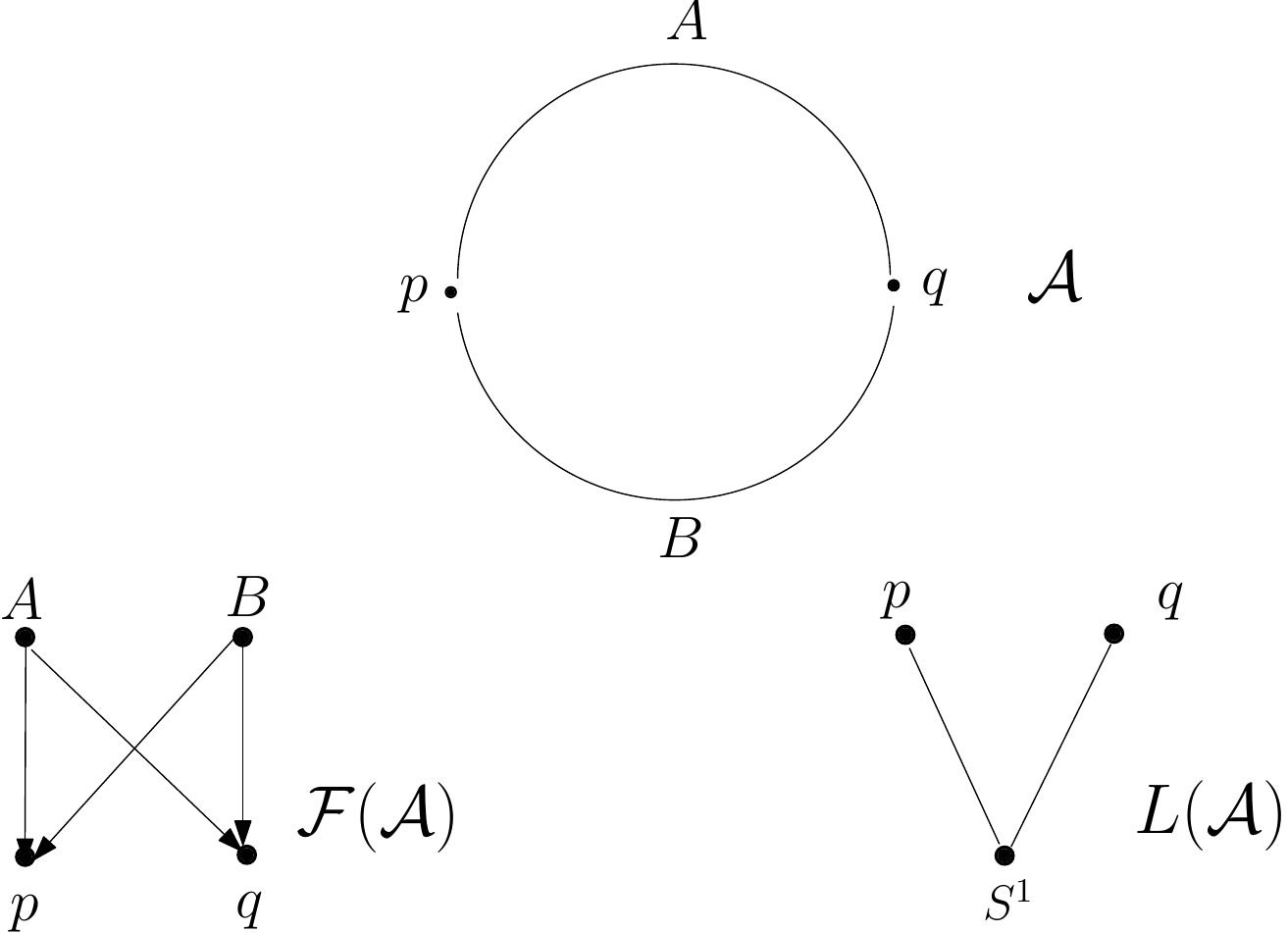} \end{center}
    \caption{Arrangement of $2$ points in a circle.}     \label{figdef1}     \end{figure}
\ee

\be{ex42}As a $2$-dimensional example consider an arrangement of $2$ great circles $N_1, N_2$ in $S^2$. Figure \ref{sphere01} shows this arrangement and the related posets. The face poset has two $0$-cells, four $1$-cells and four $2$-cells. Also note that the order complex of the face poset has the homotopy type of $S^2$. 
\begin{figure}[!ht]
\begin{center}\includegraphics[scale=0.6,clip]{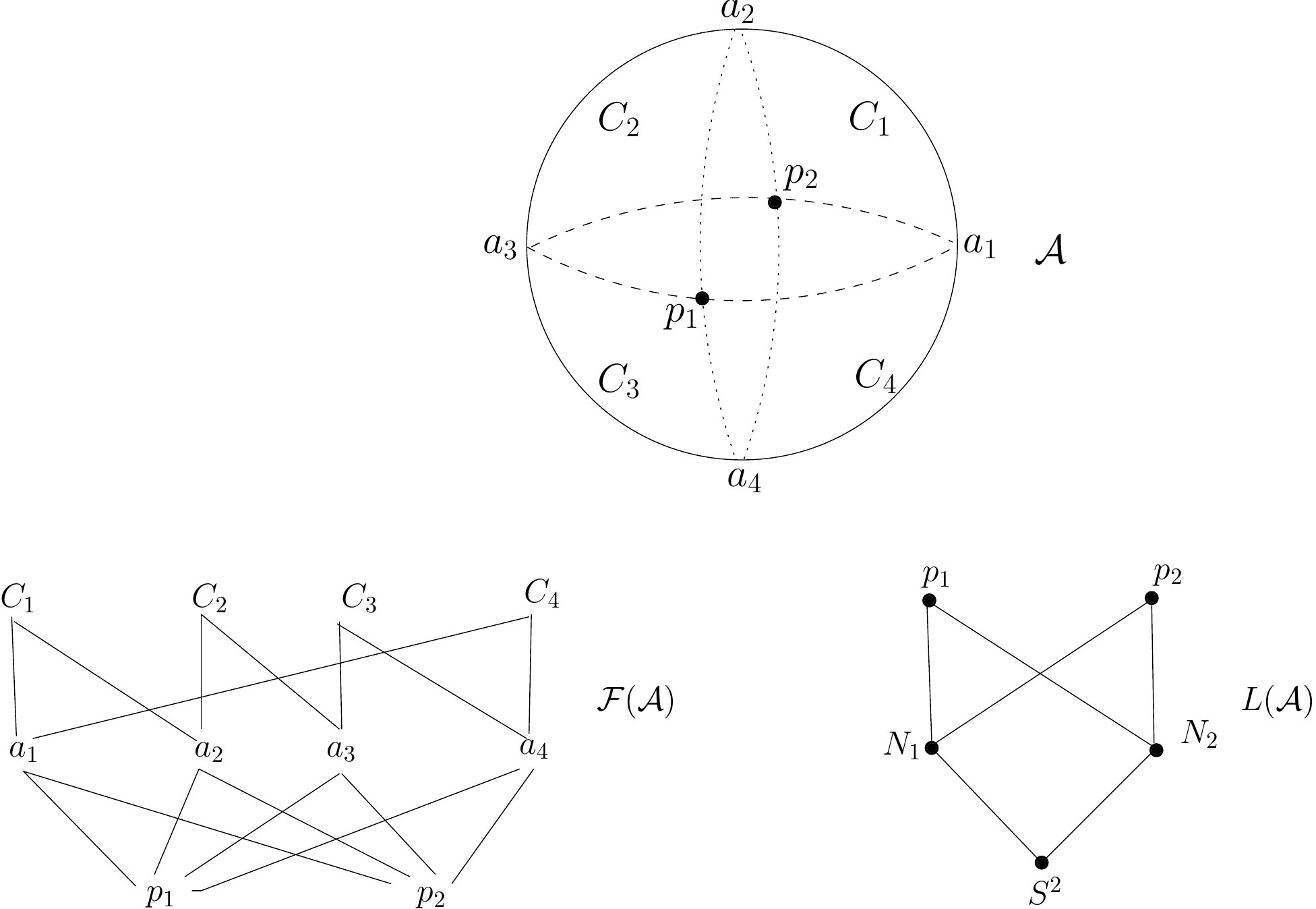}\end{center}
    \caption{Arrangement of $2$ circles in a sphere.} \label{sphere01}     \end{figure} \ee
    
\be{ex44}Now consider the $2$-torus $T^2 = \R^2/\Z^2$. A codimension-$1$ submanifold of $T^2$ is the image of a straight line (with rational slope) passing through the origin in $\R^2$ under the quotient map $\R^2\to T^2$. Let $\A$ be the arrangement obtained by projecting the lines $x = -2y,~ y = -2x,~ y =x$. These toric hyperplanes intersect in $3$ points, namely $p_1 = (0, 0), ~p_2 = (1/3, 1/3), ~p_3 = (2/3, 2/3)$. There are $6$ chambers. Figure \ref{torus} shows the arrangement ($T^2$ is considered as the quotient of the unit square) with its intersection  poset.
\begin{figure}[!ht]
\begin{center}\includegraphics[scale=0.35,clip]{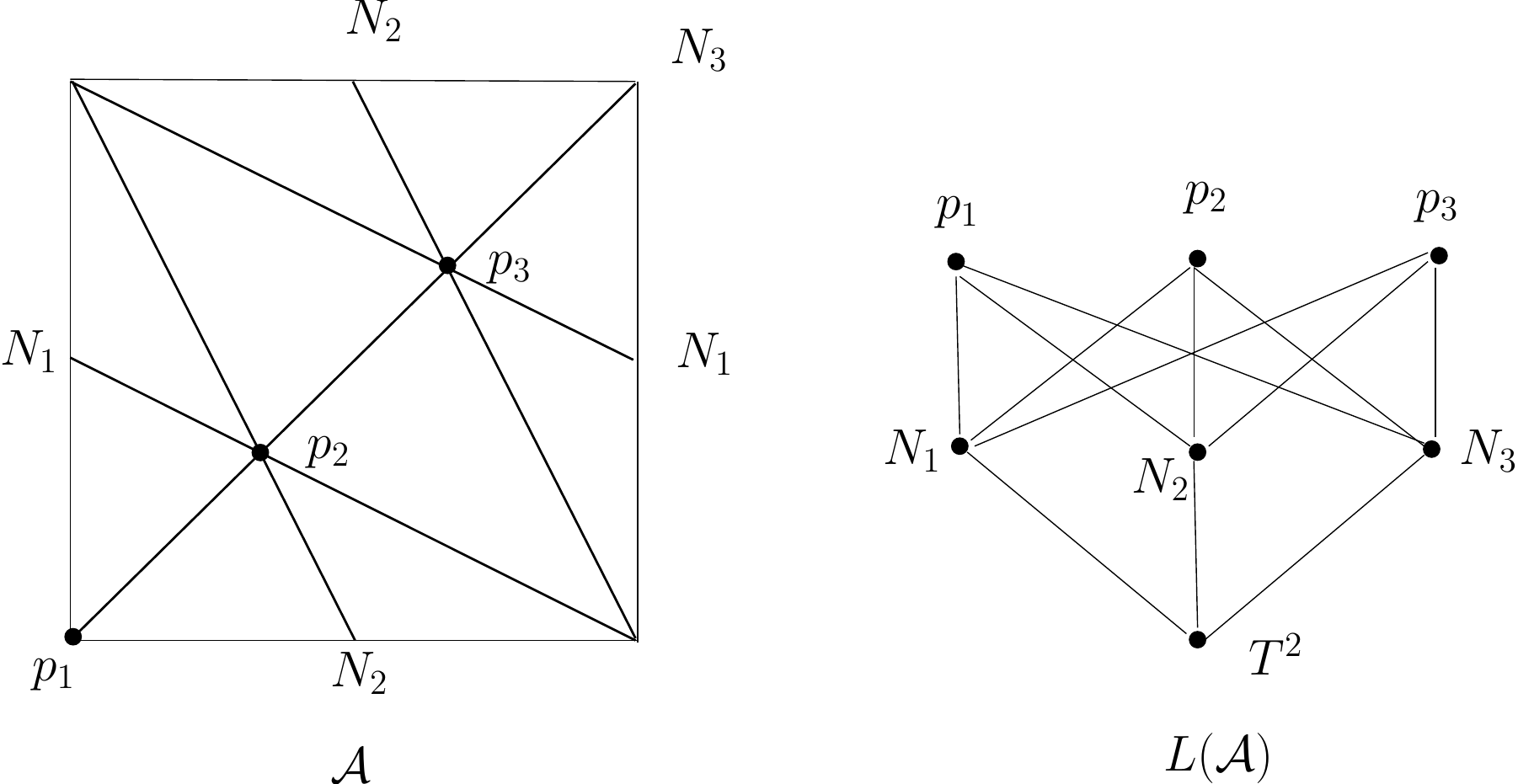}\end{center}
    \caption{Arrangement of $3$ circles in a torus.} \label{torus}     \end{figure}
\ee

Example \ref{ex44} is an example of so-called toric arrangements. This is an area which recently has seen an increase in interest. See, for example, De Concini and Processi \cite{copr05} for applications to partition functions, Ehrenborg et al. \cite{ehr09} for combinatorial perspective and also Moci \cite{moci_tutte_2009}. 

\br{remcharc} There is also a `multiplicative' way of defining toric arrangements. One can view an $l$-torus as $(S^1)^l$; i.e., a multiplicative group. Recall that the group homomorphisms from $T^l$ onto $S^1\subset \C^*$ can be described as vectors in $\Z^l$. For example, if $x_1, \dots, x_l$ are coordinates of $T^l$ and $(m_1,\cdots, m_l)\in\Z^l$ then the corresponding character is given by $(x_1, \dots, x_l)\mapsto x_1^{m_1}\cdots x_l^{m_l}$. The kernel of such a character is a codimension-$1$ submanifold. In general this kernel need not be connected; the component containing the identity is a subgroup whereas the other components are its translates. A toric arrangement is a finite collection of hypersurfaces $H_{\chi, a} := \{t\in T^l\mid \chi(t) = a \}$ where $\chi$ is a character and $a\in S^1$.  \er


\bd{def34}Let $\A$ be an arrangement of submanifolds. For a point $x\in X$ the \emph{local arrangement} at $x$ is 
\[\A_x := \{N\in \A ~|~ x\in N \}. \]
For a face $F$ the local arrangement at $F$ is
\[\A_F := \{N\in \A ~|~ F\subseteq N \}. \]
The \emph{restriction} of a local arrangement to an open set $V\subseteq X$ is
\[\A_F |_V := \{N\cap V | N\in \A_F \}. \]
\ed

\bl{lem32} Let $F$ be a face of an arrangement $\A$. Then there exists an open set $V_F\subseteq X$ containing $F$ and a map $\phi\colon V_F\to \R^l$ such that 
\begin{enumerate}
	\item $V_F\cap N = \emptyset$ for every $N\notin \A_F$;
	\item $\phi$ is a homeomorphism;
	\item $\phi$ maps $\A_F |_{V_F}$ to a central arrangement of hyperplanes in $\R^l$.
\end{enumerate} \el

\bpr Each submanifold $N\in \A$ has an open neighborhood in $X$ which is homeomorphic to $N\times (-1, 1)$ (this follows from Brown's bicollared theorem \cite[Theorem 1.7.5]{rushing73}). We call such an open neighborhood the bicollar of that submanifold. A face $F$ is homeomorphic to an open cell of some dimension and is an intersection of finitely many codimension-$1$ submanifolds. Take the open set $V_F$ to be the intersection of all the bicollars containing $F$. Moreover this $V_F$ can be chosen such that it intersects only those faces which intersect $\overline{F}$. \par
Since $X$ is a smooth manifold, there exists a homeomorphism $\phi\colon V_F\to \R^l$, say a coordinate chart. If necessary $\phi$ can be composed with a homeomorphism of $\R^l$ so that for every $N$ containing $F$, the set $N\cap V_F$ is mapped to a codimension-$1$ subspace of $\R^l$. As $\phi$ is a homeomorphism it preserves the incidence relations between the faces, thus preserving the combinatorial structure of the local arrangement $\A_F|_{V_F}$.   \epr

Intuitively every point $x\in X$ has an open neighborhood $V_x$ in $X$ which is homeomorphic to a central arrangement of hyperplanes. The hyperplanes in that arrangement correspond to $N\cap V_x$ for every $N\in \A_x$. \par 

In general $\A_F$ need not be an arrangement of submanifolds in the sense of Definition \ref{def31}. However it defines a stratification of the manifold. Let $\pi_F\colon (X, \A)\to (X, \A_F)$ be the map of the stratified spaces. It induces a map between the corresponding face categories which we again denote by $\pi_F$.

\section{The Tangent Bundle Complement} \label{TBComplement}

Recall that for a real hyperplane arrangement $\A$ the complexified complement $M(\A)$ is the complement of the union of complexified hyperplanes inside the complexified ambient vector space (Definition \ref{def4}). If one were to forget the complex structure on $\C^l$ then, topologically, it is just the tangent bundle of $\R^l$. The same is true for a hyperplane $H$ and its complexification $H_{\C}$. Hence the complexified complement of a hyperplane arrangement can also be considered as a complement inside the tangent bundle. We use this topological viewpoint to define its generalization to submanifold arrangements.

\bd{def35} Let $X$ be an $l$-dimensional manifold and $\A =\{N_1,\dots, N_k \}$ be an arrangement of submanifolds. Let $TX$ denote the tangent bundle of $X$ and let $T\A := \bigcup_{i=1}^k TN_i$. The \textbf{tangent bundle complement} of the arrangement $\A$ is defined as \[M(\A) := TX\setminus T\A. \] 
\ed

\br{remtoric}
At this point we would like to revisit the totally normal CS-space assumption in Definition \ref{def31}. We start with a familiar class of examples. Suppose that an arrangement $\A$ of hyperplanes in $\R^l$ is not essential which clearly means that the induced stratification is not totally normal (however it is cellular). If the rank of $\A$ is $k < l$ then there exists a subspace $W$ of dimension $k$ (the subspace spanned by the normals to the hyperplanes in $\A$). Projecting onto $W$ we get a new hyperplane arrangement in $W$ which we denote by $\A_W$. The arrangement $\A_W$ is the essentialization of $\A$. It is easy to see that the intersection and face posets of $\A$ and $\A_W$ respectively are isomorphic. Moreover the projection map, extended to $M(\A)$, gives a trivial vector bundle and hence
\[M(\A) \cong M(\A_W) \times \C^{l-k}.  \]
Thus it is always assumed that the hyperplane arrangements are essential. \par
Such a process can be imitated in some important examples of arrangements. We explain one such class: the toric arrangements. First note that the tangent bundle of an $l$-torus $T^l$ is homeomorphic to $(\C^*)^l$. Recall that a toric arrangement $\A$ is a collection of codimension-$1$ subtori. These subtori are the hypersurfaces defined by the vanishing of characters. The same is true for tangent bundles of the subtori: they are kernels of characters from $(\C^*)^l$ onto $\C^*$. We explain the essentialization of toric arrangements with the help of one particular example the braid arrangement.\par
Let $\A = \{ x_i x_j^{-1} = 1 ~|~ 1\leq i < j\leq l \}$ be the arrangement in $T^l$. The stratification induced by $\A$ is not totally normal and cellular. Consider the codimension-$1$ subtorus $U$ defined by 
\[ U = \{x\in T^l ~|~ x_1x_2\cdots x_l = 1 \}. \]
Let $\A_U$ denote the restriction of $\A$ to $U$. The stratification of $U$ induced by $\A_U$ is the one with the required properties. Let us see what happens to the tangent bundle complement. Define a map $\pi\colon M(\A)\to \C^*$ as $\pi(z) = z_1z_2\cdots z_l$. Check that this map is a locally trivial fiber bundle and $M(\A_U) = \pi^{-1}(1)$. This idea extends to any toric arrangement. The fiber bundle is never globally trivial; this observation corrects a statement in \cite[Remark 3.6]{dantonio_delucchi_2011}. Example \ref{ex44} above is the essentialization of the braid arrangement in $T^3$.\par 
The totally normal CS-space assumption is not very restrictive; it covers most of the important examples. Given an arbitrary collection of codimnesion-$1$ submanifolds (with locally flat intersection) it might be possible, at least in some cases, to change it to a submanifold arrangement while keeping track of the combinatorial and topological data. \er

For a point $x\in X$ let $T_x(\A) = \bigcup_{N\in \A_x} T_x(N)$.  Define the \emph{tangent space complement} at $x\in X$ as $$M(\A_x) := T_x(X)\setminus T_x(\A).$$ Then $M(\A)$ can be rewritten as 
\[M(\A) = \{(x,v)~|~ x\in X, v\in M(\A_x) \}. \]

We start the study of $M(\A)$ by first understanding the tangent space complement $M(\A_x)$. 

\bl{lem33} Let $F\in \F(\A)$ and let $\phi\colon V_F\to \R^{l}$ be a coordinate chart for an open neighborhood $V_F$ of $F$. Then for every $x\in V_F$, $$M(\A_x)\cong \phi(V_F\setminus (\A_x |_{V_F})).$$ \el

\bpr Observe that the linear map $(d\phi)_x\colon T_x(X)\to \phi(V_F)$ is an isomorphism and for the same reasons $T_x(N)\cong \phi(N\cap V_F)$ for every $N\in \A_x$. The map $\Phi\colon M(\A_x)\to \phi(V_F\setminus (\A_x |_{V_F}))$ defined by sending $(x,v)\mapsto(d\phi)_x((x, v))$ is a homeomorphism independent of the choice of coordinates.\epr

\br{rem1c3s2}
Since $\phi(\A_x |_{V_F})$ is an arrangement of hyperplanes in $\phi(V_F)$ (Lemma \ref{lem32}), $T_x(\A)$ is an arrangement of hyperplanes in $T_x(X)$. The two arrangements have isomorphic face posets as well as intersection lattices.\er

Figure \ref{local03} is a $2$-dimensional example that illustrates Lemma \ref{lem33}, showing the two homeomorphisms and chambers of the tangent space complement.

\begin{figure}[!ht]
\begin{center}
\includegraphics[scale = 0.53, clip=true]{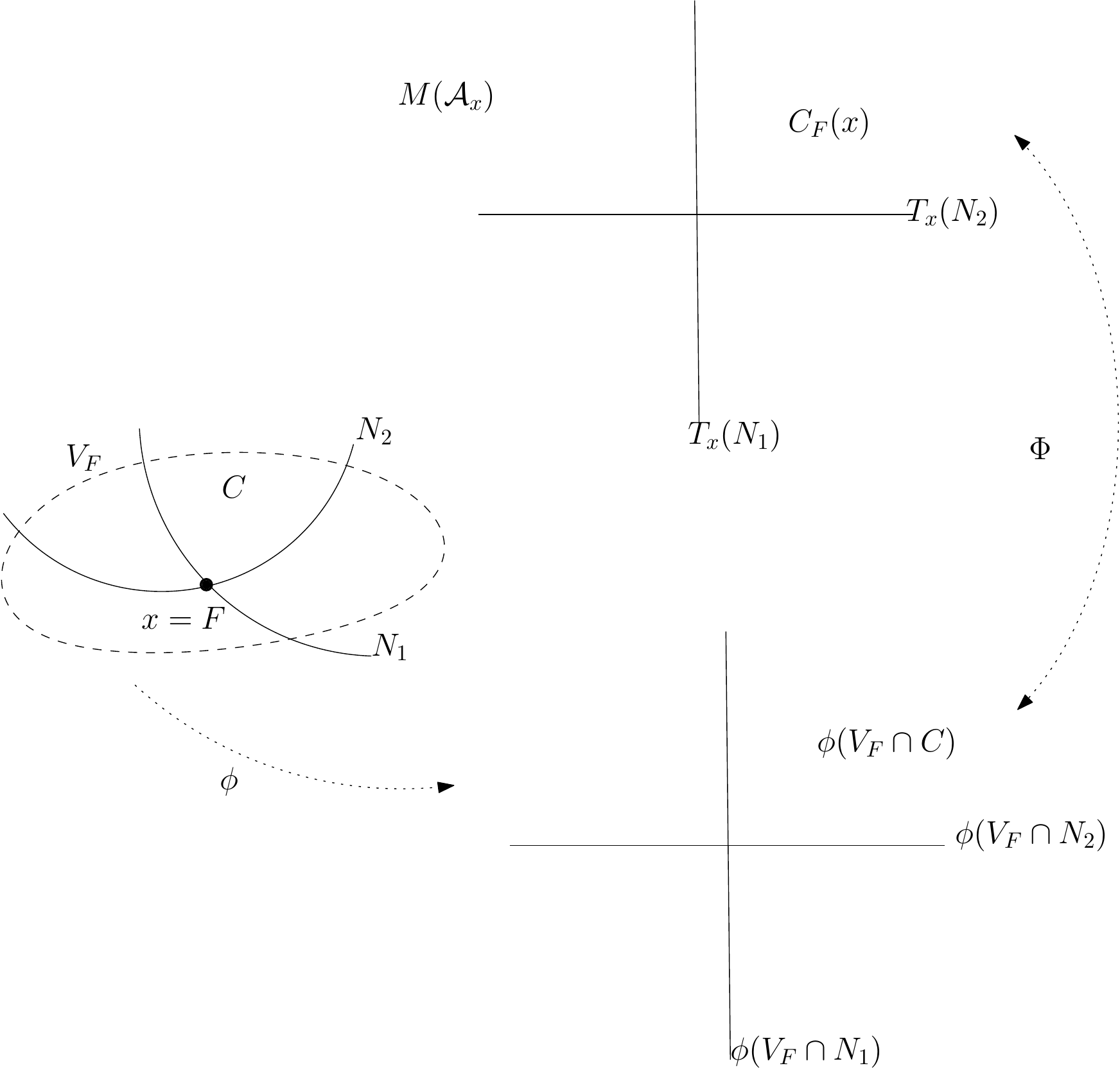}
\end{center}\caption{The tangent space complement}\label{local03}
\end{figure}

As a consequence of the above lemma chambers of the tangent space complement $M(\A_x)$ can be indexed by the chambers of $\A_x |_{V_F}$. Let us identify them. First note that if $G$ is a face such that $F\subseteq\overline{G}$ then the connected components of $V_F\cap V_G$ are in one-to-one correspondence with the morphisms from $F$ to $G$. We denote these components by $V(\alpha_{FG})$ for every morphism $\alpha_{FG} \in \F(F, G)$. Now for a chamber $C$, a morphism $\alpha_{FC}\in\F(F, C)$ and for every $x\in F$ define the following subset of $M(\A_x)$ - 
\begin{equation}U(\alpha_{FC}; x) :=\{(x,v)\in M(\A_x) ~|~ \Phi((x,v)) \in \phi(V(\alpha_{FC}))\}. \label{eq33}\end{equation}

We construct an open cover of the tangent bundle complement. The collection of open sets $\{V_F ~|~ F\in \F(A) \}$ constructed in Lemma \ref{lem32} is an open cover of $X$. These open sets can be chosen to trivialize the tangent bundle. Let $h\colon V_F\times \R^l\to \pi^{-1}(V_F)$ denote the local trivialization, where $\pi\colon TX\to X$ is the projection map. Choose a chamber $C$ such that $\F(F, C)$ is nonempty. For every point $x\in V_F$ the map $h_x := h(x, -)\colon \R^l\to \pi^{-1}(x)$ is a linear isomorphism. Let $y$ be an arbitrarily chosen point in $F$ and $U(\alpha_{FC}; y)$ be the connected component of $M(\A_y)$ corresponding to $\phi(V(\alpha_{FC}))$. Define open subsets of $\pi^{-1}(V_F)$  as
\[ W(\alpha_{FC}) := h(V_F\times h_y^{-1}(U(\alpha_{FC}; y))).\]

Note that the above subset is open and does not depend on the choice of $y$. For the sake of notational simplicity we will identify $V_F\times h_y^{-1}(U(\alpha_{FC}; y))$ with its image and will treat $W(\alpha_{FC})$ as the product $V_F\times U(\alpha_{FC}; y)$.

\bt{thm31}Let $\A$ be an arrangement of submanifolds in an $l$-manifold $X$ with $\F$ as its face category and $M(\A)$ as the associated tangent bundle complement. Then there exists a finite open cover of $M(\A)$ such that these open sets are indexed by morphisms between a face $F$ and a chamber $C$ whose closure contains $F$. Moreover, each of these open sets is contractible and so are their intersections.  \et

\bpr We show that the subsets $W(\alpha_{FC})$ constructed above satisfy the requirements. Note that by construction these are open and contractible subsets of $M(\A)$. \par

Let $(x,v)\in M(\A)$ be an arbitrary point. Suppose that $x\in F$ for some face $F$. As $v\in M(\A_x)$ we have $\phi^{-1}(\Phi(v))\in C\cap V_F$, where $C$ is some chamber whose closure contains $F$. Therefore $(x,v)\in W(\alpha_{FC})$ for some $\alpha_{FC}$. Hence the collection \[ \{ W(\alpha_{FC})~|~ \forall \alpha_{FC}\in\F(F, C)\neq\emptyset, F\in\F, C\in\Ch\}\] is an open covering of $M(\A)$.\par 
The intersection of any two such open sets is given by 
\[W(\alpha_{FC})\cap W(\alpha_{F'C'}) =  (V_F\cap V_{F'})\times ( U(\alpha_{FC}; y) \cap  U(\alpha_{F'C'}; y')).\] 
Exactly one of the following two situations can occur.
\begin{enumerate}
	\item The trivial situation being that $V_F \cap V_{F'} = \emptyset$. 
	\item $F$ is contained in the boundary of $\ol{F'}$. In that case $V_F \cap V_{F'}$ is a disjoint union of finitely many open connected subsets of $V_{F'}$; their number being the cardinality of $\F(F, F')$.
\end{enumerate}
In the second situation we have for every $y\in V_F \cap V_{F'}$ the intersection $ U(\alpha_{FC}; y)\cap  U(\alpha_{F'C'}; y)$ is a convex subset of $ U(\alpha_{F'C'}; y)$. Consequently the intersection $W(\alpha_{FC})\cap W(\alpha_{F'C'})$ is at most a disjoint union of finitely many contractible subsets. A refinement of the open cover establishes the theorem. \epr

\be{ex33}Consider an arrangement $\A = \{0\}$ in $\R$, it divides the real line in two chambers $A = (-\infty, 0)$ and $B = (0, \infty)$. The tangent bundle complement of this arrangement is the punctured plane, since $\{0\}$ is the only point whose tangent space is disconnected.
Observe that the face category of this arrangement is equivalent to the face poset.  For $p = \{0\}$ we take $V_p = (-1,1)$ and for the chambers we take $V_A = A$ and $V_B = B$ as the coordinate neighborhoods. Since there is at most one morphism between any two objects we simplify the notation for the open sets and describe them as follows: 
\begin{enumerate}
	\item $W(p,A) =\{(x,v) ~|~ x\in V_p, v < 0\}$,
	\item $W(p,B) =\{(x,v) ~|~ x\in V_p, v > 0\}$,
	\item $W(A,A) =\{(x,v) ~|~ x\in A, v\in \R\}$,
	\item $W(B,B) =\{(x,v) ~|~ x\in B, v\in \R\}$.
\end{enumerate} 
Figure \ref{wfc01} is a sketch of the corresponding open cover.
\begin{figure}[!ht] 
	\begin{center} \includegraphics[scale=0.5,clip]{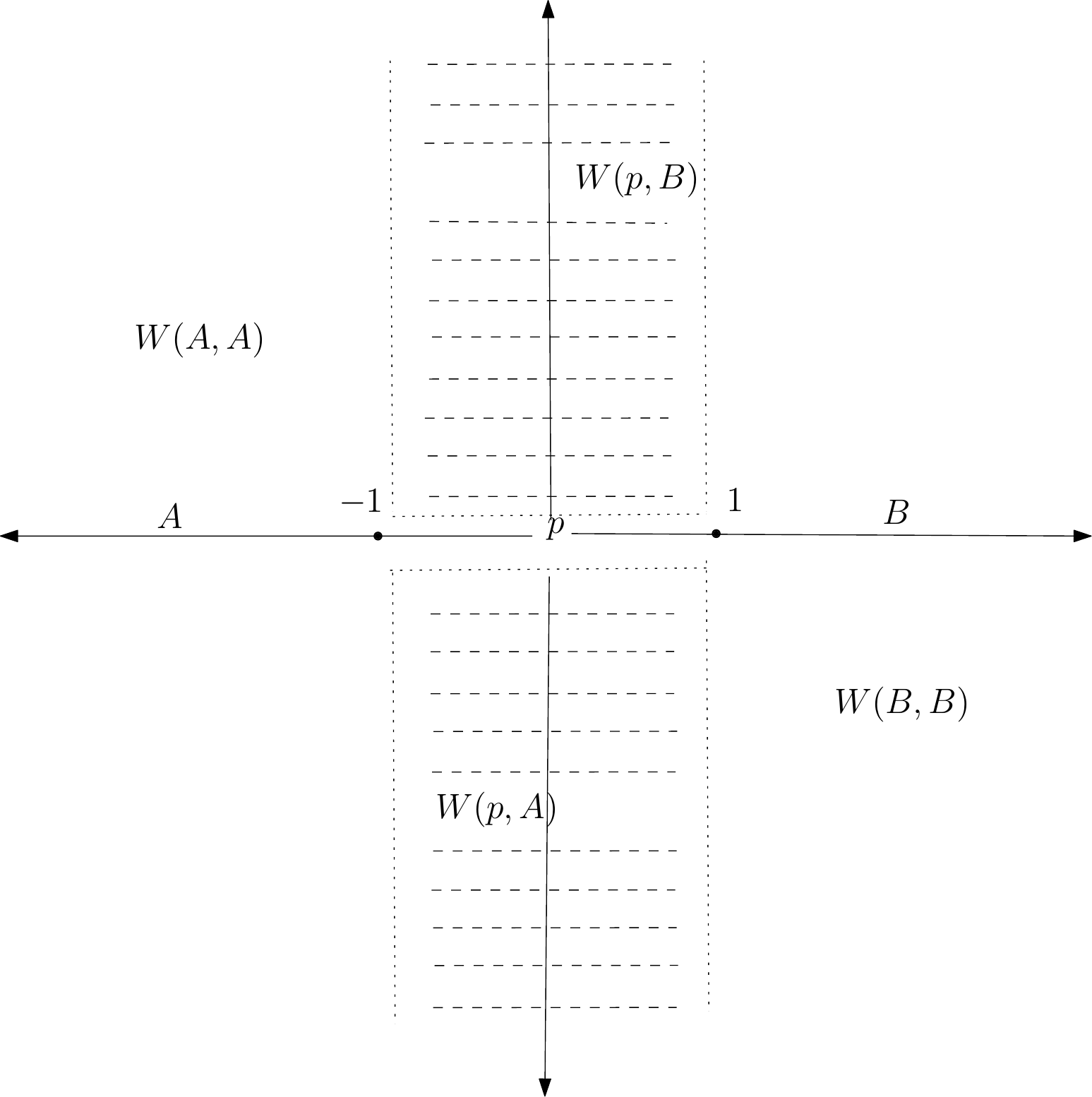}\end{center}
	\caption{Example of an open cover of the tangent bundle complement}\label{wfc01}
\end{figure} \ee

The open covering $ \{ W(\alpha_{FC})~|~ \forall \alpha_{FC}\in\F(F, C)\neq\emptyset, F\in\F, C\in\Ch\}$ constructed in Theorem \ref{thm31} satisfies the hypothesis of the Nerve Lemma. The nerve of this open cover has the homotopy type of $M(\A)$. In order to identify the simplices of this nerve we need to understand how these open sets intersect.

\bl{lem34} $W(\alpha_{FC})\cap W(\alpha_{F'C'}) \neq \emptyset$ if and only if $\F(F', F)\neq \emptyset$ and $\pi_{F}(\alpha_{FC}) = \pi_{F}(\alpha_{F'C'})$. \el

\begin{proof} From the proof of Theorem \ref{thm31} we have $V_F\cap V_{F'}\neq\emptyset$ if and only if $\F(F', F)\neq \emptyset$. We also need the other intersection in the product to be nonempty.
\begin{align*}
U(\alpha_{FC}; y)\cap  U(\alpha_{F'C'}; y) \neq \emptyset &\iff (C\cap V_F)\cap (C'\cap V_{F'})\neq\emptyset \\
&\iff (V_F\cap V_{F'})\cap (C\cap C') \neq \emptyset  \\
&\iff \pi_{F}(\alpha_{FC}) = \pi_{F}(\alpha_{F'C'}).  \qedhere
\end{align*} \end{proof}
Note that the connected components of $W(\alpha_{FC})\cap W(\alpha_{F'C'})$ can be indexed by the morphisms from $F'$ to $F$.

\begin{defn}\label{def36} Let $X$ be a smooth $l$-manifold and $\A$ be an arrangement of codimension-$1$ submanifolds. The \textbf{Salvetti category} $\mathcal{S}(\A) $ is defined as follows
\begin{enumerate}
	\item The objects of $\mathcal{S}(\A)$ are the morphisms in $\FA$ from faces to chambers, i.e., $\alpha_{FC}\in \mathrm{Mor}(\FA)$ such that $F$ is a face and $C$ is a chamber.
	\item A morphism $\beta_{F_2F_1}$ of $\FA$ defines a morphism $(\beta_{F_2F_1}, \alpha_{F_1C_1}, \alpha_{F_2C_2})$ from $\alpha_{F_1C_1}$ to $\alpha_{F_2C_2}$ in $\mathcal{S}$ if and only if $\pi_{F_1}(\alpha_{F_1C_1}) = \pi_{F_1}(\alpha_{F_2C_2})$.
	\item The composition is defined as:
	\[ (\beta_{F_3F_2}, \alpha_{F_2C_2}, \alpha_{F_3C_3})\circ (\beta_{F_2F_1}, \alpha_{F_1C_1}, \alpha_{F_2C_2}) = (\beta_{F_2F_1}\circ\beta_{F_3F_2}, \alpha_{F_1C_1}, \alpha_{F_3C_3}) \] whenever $\beta_{F_3F_2}$ and $\beta_{F_2F_1}$ are composable.
\end{enumerate}
The \textbf{Salvetti complex} $Sal(\A)$ is defined as the geometric realization of $\mathcal{S}$. \ed

\br{remDel}
Such a categorical description of the Salvetti complex first appeared in \cite{dantonio_delucchi_2011} where the authors deal with complexified toric arrangements. \er

\br{remSal}
The reader can verify that the Salvetti category is acyclic. Hence the Salvetti complex is a regular trisp. In the case when the face category is equivalent to the face poset the Salvetti category is also equivalent to a poset and the Salvetti complex is a simplicial complex.\er

\bt{them33n} Let $M(\A)$ the tangent bundle complement of an arrangement $\A$ and let $Sal(\A)$ denote the associated Salvetti complex. Then $Sal(\A) \simeq M(\A)$. \et

\bpr 
The theorem follows from the observation that the nerve of the open covering constructed in Theorem \ref{thm31} is the barycentric subdivision of the Salvetti complex.\epr
Intuitively speaking, a submanifold arrangement is constructed by gluing central hyperplane arrangements in a compatible way. It is natural to ask whether the associated tangent bundle complement is also made up by gluing local (complexified) complements. We now address this point by showing that the complement $M(\A)$ is homotopy equivalent to the homotopy colimit of a diagram of spaces over the face category. For the local arrangement $\A_F|V_F$ let $\s(\A_F)$ denote the associated Salvetti category and let $Sal(\A_F)$ denote its geometric realization. 
\bp{prop31}
With the notation as before, if $F\subseteq \ol{F'}$ are two faces then for every $\alpha_{FF'}\in\F(F, F')$ there is an inclusion $Sal(\A_{F'})\hookrightarrow Sal(\A_{F})$ as trisps. \ep

\bpr  It suffices to prove that there is an injective functor from $\s(\A_{F'})$ into $\s(\A_F)$. This follows if the face category of $\A_{F'}|_{V_{F'}}$ injects into the face category of $\A_F|_{V_F}$. But this is clear since $\F(\A){\downarrow F'}\hookrightarrow \F(\A){\downarrow F}$ and for any face $G$ the face category $\A_G|_{V_G}$ is equivalent to the comma category $\F(\A){\downarrow G}$.  \epr

\bt{thm32} Consider the diagram of spaces $\D\colon \FA^{\mathrm{op}} \to \mathrm{Top}$, given by setting $\D(F) = Sal(\A_F)$. The morphisms $\D(\alpha_{FF'})$ are the corresponding inclusions $Sal(\A_{F'})\hookrightarrow Sal(\A_F)$. Then $\mathrm{hocolim} \D \simeq Sal(\A)$. \et

\bpr Let $\E$ be a diagram of acyclic categories indexed over $\F(\A)$ such that it assigns $\s(\A_F)$ to every $F$. It follows from Lemma \ref{lem35} that $\mathrm{hocolim} \D \simeq \Delta(\mathrm{Aclim}\E)$. Hence we need to show that the functor $\phi\colon \mathrm{Plim}\E\to \mathcal{S}$ defined by $(F, (F', C))\mapsto (F',C)$ induces a homotopy equivalence $\Delta(\mathrm{Aclim} \E) \simeq Sal(\A)$, which follows from Quillen's Theorem A (Theorem \ref{qta}). \epr

\subsection{Salvetti CW-complex}\label{salcw}
We now construct a cell complex, in the spirit of Salvetti's construction, that has the homotopy type of the tangent bundle complement.\par 

Let $\A$ be an arrangement of submanifolds in an $l$-manifold $X$ and $\F(\A)$ denote the face category of the $\A$. We denote by $(X, \F(\A))$ the totally normal CS-structure of $X$ induced by $\A$. We denote by $(X, \F^*(\A))$ the \emph{dual cell structure} (see Remark \ref{rem_dual}). Every $k$-cell $F$ in $(X, \F(\A))$ corresponds to an $(l-k)$-cell $F^*$ in $(X, \F^*(\A))$ for $0\leq k\leq l$. Here $\F^*(\A)$ denotes the dual (i.e., opposite) face category. Recall that in a dual category the directions of the morphisms are reversed, i.e., $\alpha^*\in \F^*(G^*, F^*) \iff \alpha \in \F(F, G)$.  Finally, let $\prec$ denote the ordering on the poset underlying $\F^*$.\par

For the sake of notational simplicity we will denote the dual cell complex by $\F^*(\A)$ (and by $\F^*$ if the context is clear). The symbol $F^*$ will denote a $k$-cell dual to the codimension $k$-face $F$ of $\A$. Note that a $0$-cell $C^*$ is a vertex of a $k$-cell $F^*$ in $\F^*$ if and only if the closure $\overline{C}$ of the corresponding chamber contains the $(l-k)$-face $F$. \par 

Given a submanifold arrangement $\A$ in $X$ construct a totally normal CS-space $Sal(\A)$ of dimension $l$ as follows:\par

\noindent  The $0$-cells of $Sal(\A)$ correspond to $0$-cells of $\F^*$, which we denote by $\left<1; C^*, C^*\right>$.\\ 
For each $1$-cell $F^* \in \F^*$ and for every morphism $\alpha^*_{C^* F^*}$ take a homeomorphic copy of $F^*$ denoted by $\left<\alpha_{FC}; F^*, C^*\right>$. Attach these $1$-cells in $Sal(\A)_0$ (the $0$-skeleton) such that 
\[ \partial \left<\alpha_{FC}; F^*, C^*\right> = \{\left<1; D^*, D^*\right> ~|~ \F^*(D^*, F^*)\neq\emptyset \}.\] 
We put an orientation on the $1$-skeleton $Sal(\A)_1$ by directing each $1$-cell $\left<\alpha_{FC}; F^*, C^*\right>$ such that the initial vertex is $\left<1; C^*, C^*\right>$. \par

By induction, assume that we have constructed the $(k-1)$-skeleton of $S(\A)$, $1\leq k-1 < l$. To each $k$-cell $G^* \in \F^*$ and to every morphism $\alpha^*_{C^* G^*}$ ending at $G^*$ assign a $k$-cell $\left<\alpha_{GC}; G^*, C^*\right>$ that is isomorphic to $G^*$. Let $\phi(C^* G^*)\colon \left<\alpha_{GC}; G^*, C^*\right>\to Sal(\A)_{k-1}$ be the same characteristic map that identifies a $(k-1)$-cell $H^*\subseteq \partial G^*$ with the $k$-cell $\left<\beta_{HD}; H^*, D^*\right>\subseteq \partial \left<\alpha_{GC}; G^*, C^*\right>$ where $\pi_G(\alpha) = \pi_G(\beta)$. Extend the map $\phi(G^*, C^*)$ to the whole of $\left<\alpha_{GC}; G^* C^*\right>$ and use it as the attaching map, hence obtaining the $k$-skeleton. The boundary of every $k$-cell is given by 
\begin{equation}
\partial \left<\alpha_{FC}; F^*, C^*\right> = \bigcup_{\substack{G^* \prec F^*\\ \pi_G(\alpha) = \pi_G(\beta)}} \left<\beta_{GD}; G^*, D^*\right>. \label{eq1s3c3}
\end{equation}
Now we state the theorem that justifies the construction of this cell complex. 
\bt{thm1c3s3} The totally normal CS-space $Sal(\A)$ constructed above has the homotopy type of the tangent bundle complement $M(\A)$. \et
\bpr Define a functor $\Xi$ from the face category of $Sal(\A)$ to the Salvetti category as follows. $\Xi$ sends a cell $\left<\alpha_{FC}; F^*, C^*\right>$ to the object $\alpha_{FC}$. Note that, for a cell $\left<\beta_{GD}; G^*, D^*\right>$ contained in the boundary of $\left<\alpha_{FC}; F^*, C^*\right>$ the attaching maps are indexed by the morphisms in $\F^*(G^*, F^*)$. Hence $\Xi$ takes such a morphism to $(\gamma_{FG}, \alpha_{FC}, \beta_{GD})$. A routine diagram chase will convince the reader that the functor $\Xi$ is an equivalence of categories and hence induces a homeomorphism between their geometric realizations. \epr

\be{exs1p1sal} Consider the arrangement of one point in $S^1$ described in Example \ref{example1s2}. The tangent bundle complement is a once-punctured cylinder. Figure \ref{s1p1sal} shows the trisp and cellular versions of the associated Salvetti complex. 
\begin{figure}[!ht] 
	\begin{center} \includegraphics[scale=0.5,clip]{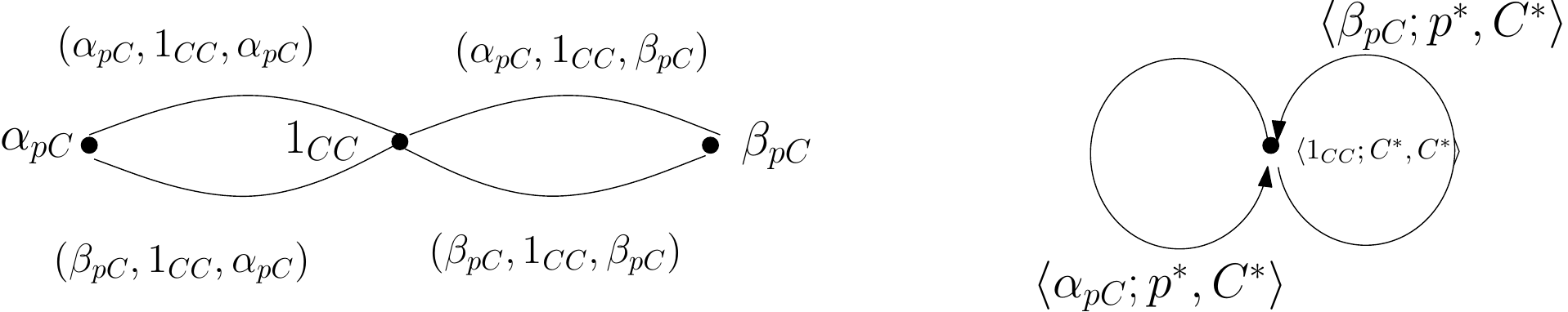}\end{center}
	\caption{Salvetti complex with the regular trisp structure and the cell structure}\label{s1p1sal}
\end{figure}
\ee

\section{Combinatorially Special Arrangements} \label{sec:csarr}

In this section we focus on a rather special class of submanifold arrangements. These arrangements are a bit closer to hyperplane arrangements in the following sense. From now on a submanifold arrangement $\A$ in $X$ will be assumed to have two additional conditions - 

\begin{enumerate}
 \item for every $N\in\A$, $X\setminus N$ is disconnected;
 \item the stratification induced by the intersections of the submanifolds imopses regular CW structure on $X$. 
\end{enumerate}

In view of assumption (1) one can consider simply connected manifolds. The regularity assumption implies that the face category of an arrangement is equivalent to the underlying poset. This assumption will help simplify the description of the Salvetti complex. Since there will be at most one morphism between any two faces we can safely drop the reference to these morphisms from the descriptions of the open covering, the Salvetti category and the Salvetti complex. Note that as a consequence the Salvetti category will be equivalent to its underlying poset and we can just refer it to as the Salvetti poset. 

\subsection{Combinatorics of the cell structure}
We say that a submanifold $N$ of $\A$ \emph{separates} two chambers $C$ and $D$ if they are contained in the distinct connected components of $X\setminus N$. For two chambers $C, D$ the set of all the submanifolds that separate these two chambers is denoted by $R(C, D)$. The \emph{distance between two chambers} is the cardinality of the set $R(C, D)$ and is denoted by $d(C, D)$. The following lemma is evident.

\bl{lem1s2c3}
Let $X$ be an $l$-manifold and $\A$ be an arrangement of submanifolds, let $C_1, C_2, C_3$ be three chambers of this arrangement. Then, 
\[R(C_1, C_3) = [R(C_1, C_2)\setminus R(C_2, C_3)] \cup [R(C_2, C_3)\setminus R(C_2, C_1)]. \] \el

Define an action of a face $F$ on a chamber $C$ as follows: 

\bd{def351} A face $F$ acts on a chamber $C$ to produce another chamber $F\circ C$ satisfying:
\begin{enumerate}
	\item $F \subseteq \overline{F\circ C}$,
	\item $\pi_F(C) = \pi_F(F\circ C)$,
	\item $d(C, F\circ C) = \hbox{~min~} \{d(C, C') ~|~ C'\in \Ch, F\subset\overline{C'} \}$. 
\end{enumerate} \ed

\bl{lem2s2c3} The chamber $F\circ C$ always exists and is unique. \el

As a consequence of the regularity assumption the open sets considered in Theorem \ref{thm31} can be written as $W(F, C)$. Lemma \ref{lem34} can now be stated as follows:

\bl{lemma34re}
$W(F, C)\cap W(F', C')\neq\emptyset$ if and only if $F'\leq F$ and $C' = F\circ C$. \el

Consequently the Salvetti poset can be described as 
\[\mathcal{S}(\A) = \{(F,C)\in \FA\times \Ch ~|~ F\leq C \}. \]

The partial order in this poset is specified in the following lemma. 

\bl{lem34n} The relation $(F_2, C_2) \leq_s (F_1, C_1)$ if and only if $F_1 \leq F_2$ and $F_2\circ C_1 = C_2$ defines a partial order on $\mathcal{S}(\A)$.\el

\begin{proof}
The arguments are similar to the proof of \cite[Lemma 3.1]{pa2}. It is obvious that the relation is reflexive and symmetric, let us check the transitivity. Pick $3$ elements such that 
$(F_3, C_3) \leq_s (F_2, C_2)$ and $(F_2, C_2) \leq_s (F_1, C_1)$. \par
The first inequality implies that $F_2 \leq F_3$ and $F_3\circ C_2 = C_3$. Similarly from the second inequality we have, $F_1 \leq F_2$ and $F_2\circ C_1= C_2$. Since $\F$ is a poset, $F_3 \leq F_1$ and $C_3 = F_3\circ(F_2\circ C_1) = F_3\circ C_1 $ which concludes the transitivity. \end{proof}

Following Theorem \ref{thm31} it is clear that the geometric realization of the Salvetti poset has the homotopy type of the tangent bundle complement. We now characterize the chains in the Salvetti poset. 

\bl{lem350}Let $\A$ be an arrangement of submanifolds, then for a chain in $\mathcal{S}$ there corresponds a chamber and a chain in $\FA$. \el 

\bpr Let $F_0\leq \cdots \leq F_k$ be a chain in $\F$ and let $C$ be a chamber. Then the Lemma \ref{lem34n} implies that $(F_k, F_k\circ C)\leq_s \cdots \leq_s (F_0, F_0\circ C)$ is a chain in $\mathcal{S}$. Moreover, using the same lemma it can be shown that every chain in $\mathcal{S}$ is of this form. \epr

\br{salsimplex}A $k$-simplex in this geometric realization is a $k$-chain in $(\mathcal{S}, \leq_s)$. Let $F_0\leq\cdots\leq F_k$ be a chain in $\FA$ and let $C$ be a chamber (such that $F_k\leq C$) then both of them determine a simplex given by
\[ (F_k,C)\leq_s\cdots\leq_s (F_0,C). \]
In fact every simplex of the nerve is of this form. \er

\be{ex4} Let $\A$ be arrangement of $2$ points in a circle. The complement $M(\A)$ in this case is $TS^1\setminus \{p,q\}$ (tangent bundle of a point is the point itself). As $TS^1\cong S^1\times\R$, the space $M(\A)$ is an infinite cylinder with $2$ punctures. Hence $M(\A)\simeq S^1\vee S^1\vee S^1$. Figure \ref{sal01} shows the geometric realization of the Salvetti poset of this arrangement.
\begin{figure}[!ht] 
\begin{center}\includegraphics[scale=0.6,clip]{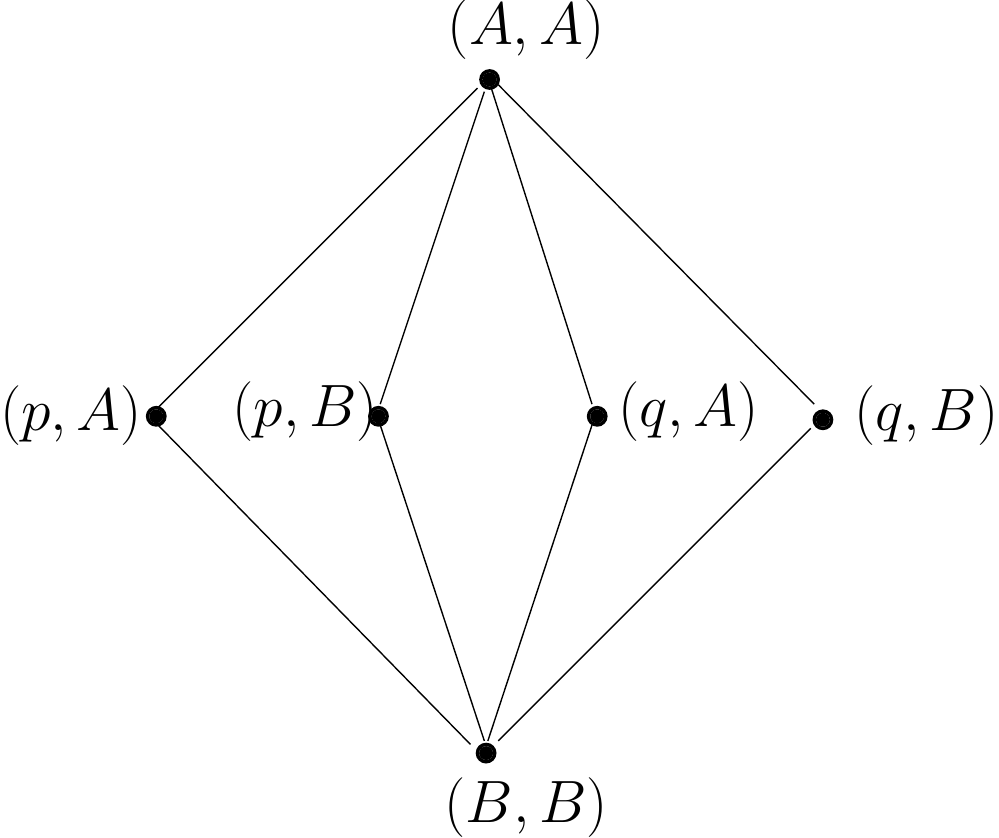}\caption{Salvetti complex}\label{sal01}\end{center}\end{figure} \ee

The cellular description of the Salvetti complex is also easier. First note that that it is now a regular CW-complex. The $k$-cells are in one-to-one correspondence with the pairs $(F, C)$ such that $F\leq C$ and the codimension of $F$ is $k$. We label this particular $k$-cell by $\left<F^*, C^*\right>$ which is homeomorphic to the dual cell $F^*$. The boundary of this cell is as follows
\begin{equation}
\partial \left<F^*, C^*\right> = \bigcup_{G^* \prec F^*} \left<G^*, G^*\circ C^*\right>. \label{eq2s3c3}
\end{equation}
We put an orientation on the $1$-skeleton $Sal(\A)_1$ by directing each $1$-cell $\left<F^*, C^*\right>$ such that the initial vertex is $\left<C^*, C^*\right>$. \par

\br{rem1s3c3}
In order to simplify the notations and the arguments we will not distinguish between the simplicial and the cellular versions of the Salvetti complex. We use the notation $[F, C]$ to denote either a cell or its barycenter and hope that the context will make it clear. Henceforth we assume that the tangent bundle complement is equipped with the above defined cell structure.\er
\be{ex1s3c3}As an example of this construction consider the arrangement of $2$ points in a circle (Example \ref{ex41}). The Figure \ref{circsal} below illustrates the dual cell structure induced by the arrangement and the associated Salvetti complex.  
\begin{figure}[!ht]
\begin{center}
	\includegraphics[scale=0.65,clip]{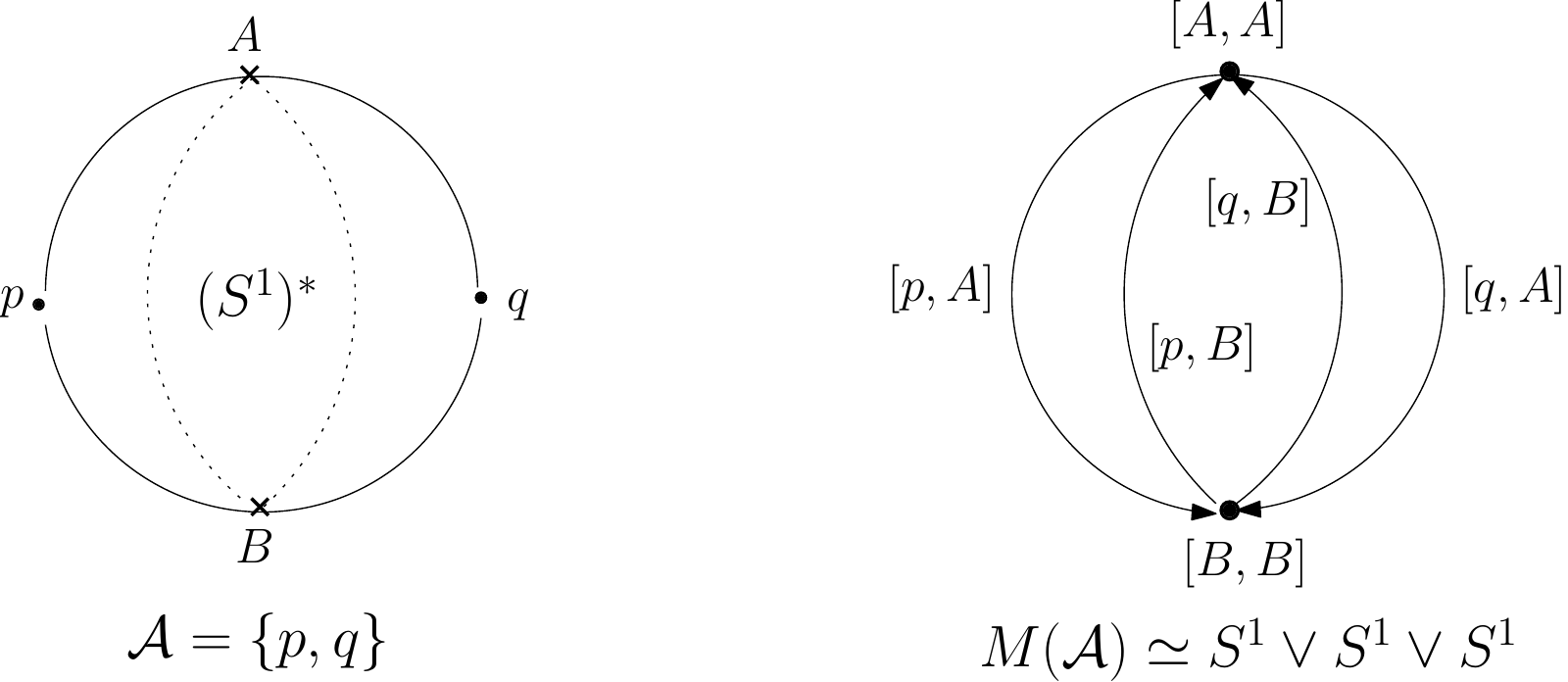}
\end{center}\caption{Arrangement in $S^1$ and the associated Salvetti complex}\label{circsal} 
\end{figure}\ee

In what follows we explore the connection between the combinatorics of the cell structure and the topology of the complement. We denote by $\F^*(\A)$ the dual face poset as well as its geometric realization. The following theorem states the connection between the manifold $X$ and the associated tangent bundle complement. 

\bt{thm2s3c3} Let $\A$ be an arrangement of submanifolds in a $l$-manifold $X$ and let $Sal(\A)$ denote the associated Salvetti complex. 
\begin{enumerate}
	\item There is a natural cellular map $\psi\colon Sal(\A)\to \F^*(\A)$ given by $[F, C]\mapsto F$. The restriction of $\psi$ to the $0$-skeleton is a bijection and in general 
	\[ \psi^{-1}(F) = \{C\in \Ch \mid C \prec F \}.\] 
	\item  For every chamber $C$ there is a cellular map $\iota_C\colon \F^*(\A)\to Sal(\A)$ taking $F$ to $[F, F\circ C]$ which is an embedding of $\F^*(\A)$ into $Sal(\A)$, and \[Sal(\A) = \bigcup_{C\in \Ch} \iota_C(\F^*).\] 
	\item The absolute value of the Euler characteristic of $M(\A)$ is the number of bounded chambers of $\A$.
\item Let $T\A$ denote the union of the tangent bundles of the submanifolds in $\A$ then,
\[\mathrm{rank}~\tilde{H}^i(TX, T\A) = \begin{cases} |\chi(M(\A))| &\hbox{~\emph{if}~} i = l, \\ 0 &\hbox{~\emph{otherwise.}}\end{cases} \]
\end{enumerate} \et

\begin{proof}
(1) and (2) follow from the simple observation that $\psi\circ\iota_C$ is identity on $\F^*(\A)$. This also implies that $X$ is homeomorphic to a retract of $M(\A)$.\par  
We prove (3) by explicitly counting cells in the Salvetti complex. The Euler characteristic of a CW complex $K$ is equal to the alternating sum of number of cells of each dimension. Given a $k$-dimensional dual cell $F$ there are as many as $|\{C\in\Ch \mid F\leq C \}|$ $k$-cells in $\psi^{-1}(F)$. Hence for a vertex $[C, C]\in Sal(\A)$ the number of $k$-cells that have $[C, C]$ as a vertex is equal to the number of $k$-faces contained $\ol{C}$. The alternating sum of numbers of such cells is $1 - \chi(Lk(C))$, where $Lk(C)$ is the link of $C$ in $\F^*(\A)$. Applying this we get 
\[\chi(Sal(\A)) = \sum_{C\in \Ch}(1 - \chi(Lk(C))). \] 
If a chamber is unbounded then $Lk(C)\simeq \mathbb{B}^{l-1}$ and on the other hand if it is bounded then $Lk(C)\simeq S^{l-1}$. Hence we have, 
\begin{align*}
		\chi(Sal(\A)) &= \sum_{C\in \Ch}(1 - \chi(Lk(C)))  \\
		              &= \sum_{C {~\mathrm{unbounded}}}(1 - 1) + \sum_{C {~\mathrm{bounded}}}(1 - [1 + (-1)^{l-1}]) \\
					  &= (-1)^l \sum_{C {~\mathrm{bounded}}}1.
\end{align*}
Hence, \[ \chi(M(\A)) = (-1)^l \hbox{(number of bounded chambers)}. \]
Let $\bigcup\A$ denote the union of submanifolds in $\A$. Since $\A$ induces a regular cell decomposition it has the homotopy type of wedge of $(l-1)$-dimensional spheres. The number of spheres is equal to the number of bounded chambers. Claim (4) now follows from the homeomorphism of pairs $(TX, T\A) \cong (X, \bigcup\A )$.
\end{proof} 

\subsection{Metrical-hemisphere complexes} \label{mhcSection}

We now take a closer look at the cells of a Salvetti complex. Our aim is to understand how the combinatorial properties of the Salvetti complex associated to an arrangement of hyperplanes generalize in the context of submanifold arrangements. In particular, we show that the combinatorics is similar to that of a zonotope. Recall that a central arrangement of hyperplanes decomposes the ambient Euclidean space into open polyhedral cones. As a matter of fact every hyperplane arrangement is a normal fan of a very special polytope known as the zonotope. Zonotopes can be defined in various ways: for example, projections of cubes, Minkowski sums of line segments, dual (polar) of hyperplane arrangements etc. For more on the relationship between zonotopes and hyperplane arrangements see Ziegler \cite[Lecture 7]{zig95}.

\bd{defc1} A \emph{zonotope} is a polytope all of whose faces are centrally symmetric (equivalently every $2$-face is centrally symmetric). A \emph{zonotopal cell} is a closed $k$-cell such that its face poset is isomorphic to the face poset of a $k$-zonotope. \ed

The face poset of a zonotope has some special combinatorial properties, the most important of which is the product structure. This product is basically the one on the face poset of a hyperplane arrangement or on the set of covectors of an oriented matroid. We show that the dual cell structure of a submanifold arrangement and the cell structure of its associated Salvetti complex enjoy similar combinatorial structure. We do not use the language of zonotopes to explain this combinatorics (since it is exclusive to Euclidean settings); instead we use the language of \emph{metrical-hemisphere complexes}. These cell complexes generalize zonotopes from a topological viewpoint. The metrical-hemisphere complexes (MH-complexes for short) were first introduced by Salvetti in \cite{salvetti93} where he generalized his construction and proved an analogue of Deligne's theorem for oriented matroids. \par

Let $Q$ denote a connected, regular, CW complex (and $|Q|$ be the underlying space). The $1$-skeleton of such a complex $Q$ is a graph $G(Q)$ with no loops (abbreviated to $G$ if the context is clear). The vertex set of this graph will be denoted by $VG$ and the edge set by $EG$. An edge-path in $G(Q)$ is a sequence $\alpha = (l_1, \dots, l_n)$ of edges that correspond to a connected path in $|Q|$. The inverse of a path is again a path $\alpha^{-1} = (l_n,\dots, l_1)$. Two paths are composed by concatenation if the ending vertex of one of the paths is the starting vertex of another. The distance $d(v, v')$ between two vertices will be the least of the lengths of paths joining $v$ to $v'$. Given an $i$-cell $e^i\in Q$, let $Q(e^i) := \{e^j\in Q : |e^j| \subset |e^i| \}$ and $V(e^i) = VG\cap Q(e^i)$.

\bd{def2s3c3}A regular CW complex $Q$ is a QMH-complex (quasi-metrical-hemisphere complex) if there exist two maps 
\[\mino, \maxo \colon VG\times Q\to VG \] 
such that for all $v\in VG, e^i\in Q$ following properties are satisfied.
\begin{enumerate}
	\item $\mino(v, e^i)\in V(e^i)$ and $d(v, \mino(v, e^i)) = \hbox{minimum} \{d(v, u) | u \in V(e^i) \}$.
	\item $\maxo(v, e^i)\in V(e^i)$ and $d(v, \maxo(v, e^i)) = \hbox{maximum} \{d(v, u) | u \in V(e^i) \}$.
	\item $d(v, \maxo(v, e^i)) = d(v, u) + d(u, \maxo(v, e^i))$ for all $u\in V(e^i)$.
\end{enumerate}  \ed

This definition imposes a strong restriction on the $1$-skeletons of such complexes (see \cite[Proposition 1]{salvetti93}). 

\bl{lem1s3c3} If $Q$ is a QMH-complex then each circuit in $G$ has an even number of edges. Moreover if $Q$ is homeomorphic to a $k$-ball then it is a zonotopal cell. \el
%
%
%
%
For any $e^i\in Q$, indicate by $G(e^i)\subset G(Q)$ the subgraph corresponding to the $1$-skeleton of $e^i$ and by $d_{G(e^i)}$ the distance computed using $G(e^i)$. 

\bd{def3s3c3}A regular CW complex will be called an LMH-complex (local-metrical-hemisphere complex) if each $Q(e^i)$ is a QMH-complex with respect to $d_{G(e^i)}$. Moreover, the following compatibility condition also holds: 
if $e^k\in Q(e^i)\cap Q(e^j), v\in V(e^i)\cap V(e^j)$ then 
\begin{equation*} 
\mino_{(e^j)}(v, e^k) = \mino_{(e^i)}(v, e^k), \quad \maxo_{(e^j)}(v, e^k) = \maxo_{(e^i)}(v, e^k).
\end{equation*} 
Here, $\mino_{(e^j)}, \maxo_{(e^j)}$ are defined similar to $\mino, \maxo$ but using $d_{G(e^j)}$. Finally, $Q$ is an \textbf{MH-complex} if $Q$ is both a QMH-complex and an LMH-complex and for all $e^i\in Q, e^j\in Q(e^i), v\in V(e^i)$
\begin{equation*}
\mino(v, e^j) = \mino_{(e^i)}(v, e^j), \quad \maxo(v, e^j) = \maxo_{(e^i)}(v, e^j). \end{equation*} \ed

\br{rem3s3c3}
Note that the $1$-skeleton of an MH-complex has very special properties with respect to the distance. It is not enough to have a cell complex all of whose cells are zonotopal. Here are two examples that illustrate the special nature of MH-complexes.\par 
%
Consider a $2$-dimensional cell complex made up of two $1$-cells and  two $2$-cells. The $1$-cells are attached to an octagonal $2$-cell. There is a one more trapezoidal $2$-cell whose three $1$-cells are attached to three $1$-cells in the boundary of the octagonal cell as shown in Figure \ref{nolmh}. The resulting complex is QMH but not LMH. Consider the $1$-cell labeled by $e$ in the figure. There are two vertices, namely $v_1, v_2$, in its boundary. Considering $e$ as a member of the trapezoidal cell we see that the vertex $v_1$ is closest to the vertex $v_4$. On the other hand as a member of the octagonal $2$-cell vertex $v_2$ is closest to $v_4$.
\begin{figure}[!ht] \begin{center}
	\includegraphics[scale=0.55,clip]{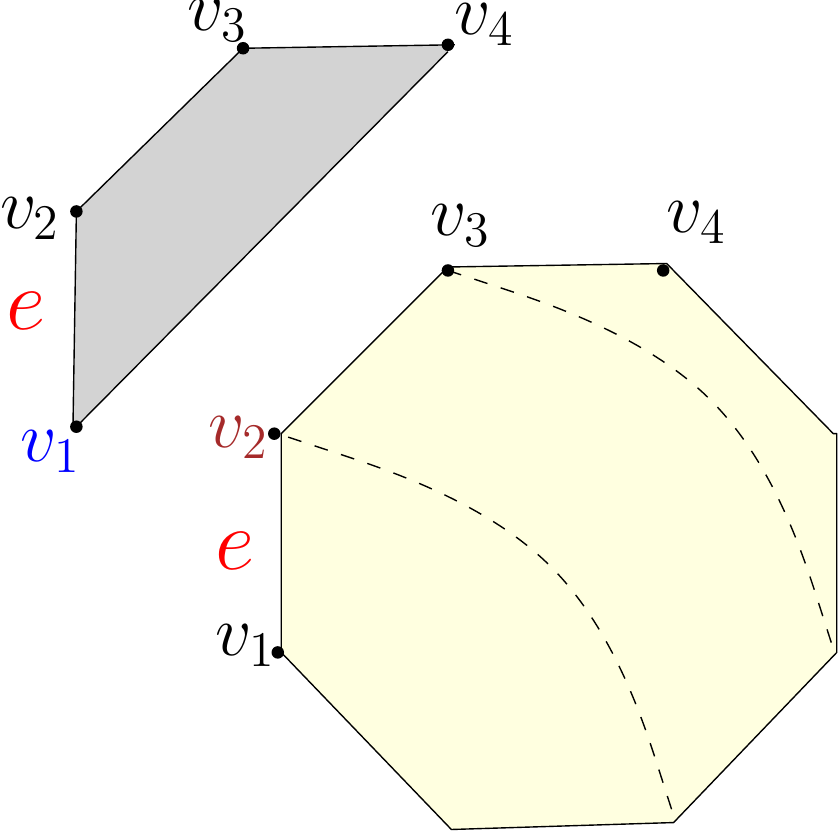}
\end{center} \caption{A QMH complex without LMH structure.} \label{nolmh}\end{figure}

The next example shows a cell complex (Figure \ref{nomh}) obtained by removing the trapezoidal $2$-cell from the first example. The resulting cell complex is both QMH and LMH but not an MH-complex. Consider the $1$-cell labeled by $e$, there are two boundary vertices $v_1, v_2$. Considering $e$ as a member of the octagonal cell the vertex $v_2$ is closest to $v_3$. But in the whole complex the boundary vertex of $e$ closest to $v_3$ is $v_1$.
\begin{figure}[!ht]\begin{center}
	\includegraphics[scale=0.5,clip]{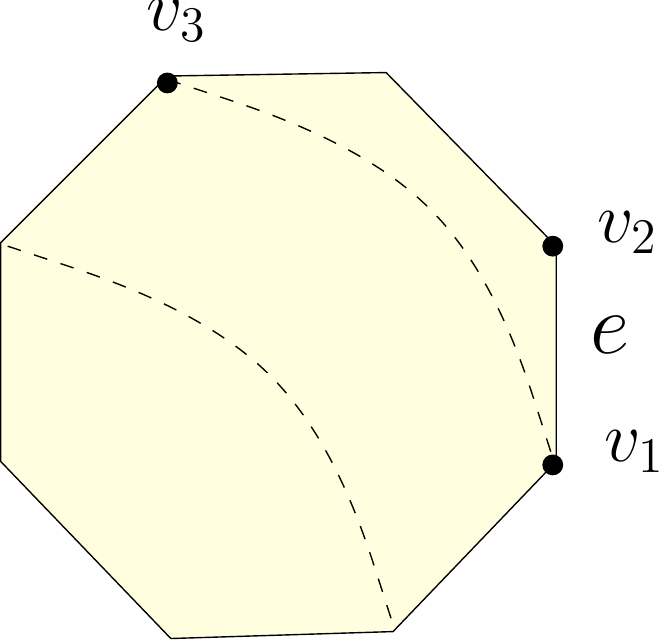}
\end{center}\caption{A QMH and LMH-complex which is not an MH-complex}\label{nomh}\end{figure}   \er

The following lemma establishes the combinatorial connection between zonotopes and MH-complexes. It states that the distance between any two vertices is the same no matter how it is measured, either locally or globally (see \cite[Proposition 5]{salvetti93}).
\bl{lem2s3c3} Let $Q$ be an MH-complex, $e^i\in Q, v,v'\in V(e^i)$. Then $d(v, v') = d_{G(e^i)}(v, v')$. \el

\bpr Let $\alpha = (l_1,\dots, l_n)$ be a minimal path of $G(e^i)$ between $v$ and $v'$ (so $d_{G(e^i)}(v, v') = n$). Let $v_{j-1}, v_j$ be the vertices of $l_j$ ordered according to the orientation of $\alpha$ from $v$ to $v'$. Since $\alpha$ is minimal in  $G(e^i)$ and $Q$ is an MH-complex we have,
\[\mino_{(e^i)}(v, l_j) = v_{j-1} = \mino(v, l_j) \]
Hence, $d_{G(e^i)}(v, v_j) = d_{G(e^i)}(v, v_{j-1}) + 1$ and $d(v, v_j) = d(v, v_{j-1}) + 1$ for $j = 1,\dots, n$.\epr

We now state the theorem that generalizes the relationship between hyperplane arrangements and zonotopes. Recall that $\F^*$ denotes the dual cell structure and that chambers of $\F$ are vertices of $\F^*$.

\bt{thm3s3c3}
Let $X$ be a smooth manifold of dimension $l$, $\A$ denote an arrangement of submanifolds. The dual cell complex $\F^*$ is an MH-complex. \et
\bpr First, we need to define the two maps $\mino, \maxo$ and then show that they are well defined. Let $F^i$ be an $i$-cell and $C$ be a vertex of $\F^*$ then,
\[\mino(C, F) := F\circ C. \]
Using the same strategy as in the proof of Lemma \ref{lem2s2c3} we can identify a unique chamber (of $\A$) whose closure contains (dual of) $F$ and is farthest from $C$, denote this chamber by $F\ast C$ and 
\[ \maxo(C, F) := F\ast C. \]
This shows that the maps $\mino, \maxo$ are well defined. Second, note that a path between two vertices $C, C'$ has minimal length (among all paths from $C$ to $C'$) if and only if it crosses the faces that separate $C$ from $C'$ exactly once and does not cross any other face. The distance between any two vertices of $\F^*$ is thus 
\[d(C, C') :=  |R(C, C')|. \]
Observe that if $F \subseteq N_1\cap\cdots\cap N_r$ then $R(F\circ C, F\ast C) = \{N_1,\dots, N_r \}$ and if $D\in V(F)$ then by Lemma \ref{lem1s2c3}
\[R(F\circ C, D)\textstyle \cup R(F\ast C, D) = R(F\circ C, F\ast C) \hbox{~and~} R(F\ast C, D)\cap R(F\circ C, D) =\emptyset. \]
Moreover, 
\[ R(C, D) = R(C, F\circ C)\textstyle\cup R(F\circ C, D),~~ R(C, F\circ C)\cap R(F\circ C, D) = \emptyset \] 
and 
\[ R(C, D) = R(C, F\ast C)\textstyle\cup R(F\ast C, D),~~ R(C, F\ast C)\cap R(F\ast C, D) = \emptyset.\]

Using the last equality we see that $\F^*$ is a QMH-complex. The other compatibility conditions also follow easily. \epr

In a nutshell the above theorem justifies the following statement. MH-complexes are to submanifold arrangements as zonotopes are to hyperplane arrangements.\par 

In case of hyperplane arrangements, each cell of the associated Salvetti complex is zonotopal (see \cite[Proposition 5.7]{bz92}). We generalize this claim and also show that in general the Salvetti complex is a MH-complex. 

\bt{thm4s3c3} Let $\A$ be an arrangement of submanifolds in a manifold $X$. The associated Salvetti complex $Sal(\A)$ is an MH-complex. \et

\bpr The $1$-skeleton of $Sal(\A)$ is obtained by `doubling' the edges in the $1$-skeleton of $\F^*$. Hence the distance between any two vertices of $Sal(\A)$ is same as the distance between the corresponding vertices of $\F^*$. Also, by construction, there is a one-to-one correspondence between vertices of $[F, C] $ and the vertices of $F$ for all $F\in \F^*$. \epr

\br{rem3news3c3}Note that in the construction of the Salvetti complex as mentioned before we used the $\mino(F, C) = F\circ C$ map. It is possible to repeat the same construction using the $\maxo$ map. However these two complexes are isomorphic. \er

As a consequence of Lemma \ref{lem1s3c3}, every $2$-cell of the Salvetti complex (of an arrangement) is combinatorially equivalent to an oriented polygon with an even number of edges. This observation together with the condition $(3)$ in Definition \ref{def2s3c3} (QMH-complex) proves the following. 

\bc{cor2s3c3}
In the Salvetti complex $Sal(\A)$, the oriented $1$-skeleton of a $k$-cell $[F, C]$ is composed of edge-paths going from $[C, C]$ to $[F\ast C, F\ast C]$. All these paths have same lengths and are directed away from $[C, C]$. \par
\ec

\begin{proof} 
Assume that there is an edge $[G, C']$ contained in the boundary of $[F, C]$. Hence, $C' = G\circ C$ and the result follows because of condition (3) defining the QMH-complex. \end{proof}

\section{Concluding Remarks}\label{endremarks}
There are several directions one can take from here. Probably the next step should be to understand the cohomology ring of the tangent bundle complement. For hyperplane arrangements there is a beautiful theory of cohomology of the complement, see Orilk and Terao \cite{orlik92} for details. However, till date, there is very little understanding of how the combinatorics of submanifold arrangements determines the cohomology ring of the complement. We would like to direct the reader to the author's thesis \cite[Section 3.7]{deshpande_thesis11} for some calculations of cohomology groups using the Bousfield-Kan spectral sequence. The examples analyzed so far suggest that there is a finer grading of cohomology groups indexed by the intersection data \cite[Conjecture 3.7.8]{deshpande_thesis11}. \par

A current work in progress is to extend the language of Salvetti-type diagram models  developed by Delucchi \cite{del1} in order to describe the connected covers of the tangent bundle complement in terms of the face category. It will be interesting to see how the MH-complex construction works in this context; it might lead to a generalization of oriented matroids. \par 

One of the motivations to study hyperplane arrangements comes from its natural connection with the reflection groups and their associated Artin groups. To every finite reflection group there corresponds an arrangement of reflecting hyperplanes. There is a fixed-point free action of the reflection group on the complexified complement. The fundamental group of the resulting orbit space is the Artin group (of finite type) associated with the (Coxeter presentation of the) reflection group. The topological properties of this space offer an insight into algebraic properties of the Artin group. Deligne's work in \cite{deli72} is pioneering in that respect. It has greatly influenced research in hyperplane arrangements and geometric group theory. For example, existence of a bi-automatic structure in finite type Artin groups as proved by Charney \cite{charney92} and the discovery of Garside groups by Dehornoy and Paris \cite{dehornoy99} are consequences of Deligne's work.\par  

The present work can be seen as a step towards generalizing Artin groups. There is a well developed theory of geometric reflection groups and more generally of manifold reflection groups, for example, see Davis \cite[Chapters 6, 10]{davisbook08}. A (separating) reflection on a smooth manifold is a diffeomorphism of order $2$ such that it fixes a codimension-$1$ submanifold which separates the ambient manifold. A discrete subgroup of diffeomorphisms generated by finitely many reflections is called the Coxeter transformation group. Such a group has a Coxeter presentation. The fixed point set is indeed a submanifold arrangement. The transformation group acts fixed point freely on the tangent bundle complement. The fundamental group of the complement is the analogue of pure Artin groups whereas the fundamental group of the orbit space is the analogue of Artin groups. It will be interesting to explore these generalized Artin groups using the results proved in this paper.

\bibliographystyle{abbrv} 
\bibliography{arrman} 

\end{document}